



\input amstex
\documentstyle{amsppt}

\magnification\magstep1

\define\p{\Bbb P}

\define\a{\Bbb A}

\redefine\c{\Bbb C}

\redefine\o{\Cal O}

\define\q{\Bbb Q}

\define\f{\Bbb F}

\define\z{\Bbb Z}

\define\qtq#1{\quad\text{#1}\quad}

\define\section#1{

\bigpagebreak{\smc#1}\bigpagebreak

}

\define\rup#1{\ulcorner #1\urcorner}

\define\rdown#1{\llcorner #1\lrcorner}

\define\demop{\demo{Proof}}

\define\spec{\operatorname{Spec}}

\define\red{\operatorname{red}}

\define\stab{\operatorname{Stab}}

\define\supp{\operatorname{Supp}}

\redefine\hom{\operatorname{Hom}}

\define\len{\operatorname{length}}

\define\im{\operatorname{im}}

\define\chr{\operatorname{char}}

\define\deq{:=}

\define\norm#1{\vert\vert#1\vert\vert}

\topmatter
\title Quotient Spaces Modulo Algebraic Groups
\endtitle
\author J\'anos Koll\'ar
\endauthor
\address Department of Mathematics, University of Utah, Salt Lake City
UT 84112, USA
\endaddress
\email kollar\@{}math.utah.edu \endemail
\endtopmatter

\head 1. Introduction
\endhead

In algebraic geometry one often encounters the problem of taking
the quotient of a scheme by a group. Despite the frequent occurrence of such
problems, there are few general results about the existence of such
quotients. These questions come up again and again in the theory of
moduli spaces. When we want to classify some type of algebraic objects, say
varieties or vector bundles, the classification usually proceeds in two steps.

First we classify  the objects with some extra structure added. In
classifying varieties  we first parametrize  embedded varieties in a fixed
projective space $\p^N$.  In studying vector bundles we first describe them
using a basis in the vector space  $V$ of global sections.

Second, we have to get rid of the effect of the extra structure. In the case
of varieties this usually means dividing out by the automorphism group of
$\p^N$. In the case of vector bundles we   take the quotient by the
automorphism group of the space of sections $V$.

There have been several attempts to give a general theory of quotients in
algebraic geometry.

Mumford's geometric invariant theory is designed to
construct quotients by group actions. This approach works very successfully
for the quotient problems arising in the theory of vector bundles. Actually,
in this area geometric invariant theory solves two problems at once. It
predicts the best equivalence relation by which a quotient might exist (a
problem with semistable objects), and then it proceeds to construct the
quotient.

In the theory of the moduli of varieties, geometric invariant theory
has been less successful. It does not seem to predict the  correct class of
varieties where a good quotient should exist. The existence of the quotient
is rather difficult to establish, and it has been done only over $\c$
\cite{Viehweg95}.

Artin developed the theory of algebraic spaces in order to tackle more general
quotient problems. It was noticed early on that the quotient of a scheme by a
finite group is, in general, not a scheme. (The simplest example is smooth
of dimension 3). Naively one can think of algebraic spaces as
quotients of schemes by finite groups (this is completely correct for normal
algebraic spaces). This approach is   inconvenient in practice, and the
working definition is different \cite{Knutson71}.

The theory of algebraic spaces is rather successful in handling general
quotient problems. For instance, the quotient of an algebraic space by a flat
equivalence relation is again an algebraic space \cite{Artin74}.  This
implies that if $g:S'\to S$ is faithfully flat then  descent data for
algebraic spaces over $S'$ are always effective.

Actions of groups do not produce flat equivalence relations in general, but
they are not far from it. Despite this  close relationship, it has
been rather difficult to prove quotient theorems modulo group actions.

One can not expect the existence of quotients for arbitrary actions.  If a
group acts on a topological space, a quotient always exists, but it may not be
separated.   The notion of a group acting ``properly" was developed
to ensure the separatedness of a quotient. The concept of a proper
action is rather easy to translate to algebraic geometry (2.3).
Unfortunately, in concrete examples it is frequently hard to check if
the group action is proper or not. One difficulty is that being proper is not
a local condition in neighborhoods of orbits.
Moreover, in considering moduli of vector bundles, the occurring group actions
are oftentimes not proper.

The situation is completely
different in the theory of polarized varieties, and especially
for canonically polarized varieties. The group actions are almost always
proper and checking properness is usually the easiest task. Therefore it
is reasonable to pay special attention to  proper group actions.

The following emerged as  a natural  problem in this area:

\proclaim{1.1 Conjecture}  Let $X$ be a scheme on which an algebraic group
scheme $G$ acts properly. Then the geometric quotient  $X/G$ exists (2.7.3).
\endproclaim

\demo{1.2 Clarification} In this generality the conjecture definitely fails.
$G$  and   $X$  must satisfy some    technical assumptions.
If $G/S$ has varying fiber dimension, then  the upper semicontinuity
of dimensions shows that a quotient cannot exist in most cases.  Therefore,
  we have to assume that  $G\to S$ is flat or
at least universally open.

 \cite{Nagata69} constructed
examples of   Artin rings in characteristic $p$  with a
$\z_p$-action whose invariant ring  is not  Noetherian (6.5.1).
(The better known constructions of \cite{Nagata59} are  not relevant
here. In those examples the quotient is an algebraic space of finite type.) To
exclude similar cases, we assume that
$G$ and $X$ are of finite type over a base scheme.

The conjecture may  hold with these assumptions.
\enddemo

\demo{1.3 Main Cases} In most  of the applications that have been
developed so far, one needs this conjecture only in case $G=PGL(n,K)$ for
a field $K$ (for moduli spaces over a field) and for $G=PGL(n,\z)$ (for
arithmetic moduli questions).

Several moduli problems have been considered  where the group in
question is not reductive. For instance, the moduli of parabolic bundles
appears naturally as  a quotient by a parabolic subgroup of $GL$
\cite{Seshadri82}. In many cases the most economical embedding of a variety is
to a weighted projective space. Then we need to take the quotient  by the
automorphisms of a weighted projective space, which is not reductive
in general.
In the theory  of singularities, many of the naturally
occurring groups are unipotent \cite{Greuel-Pfister93}.
\enddemo

\demo{1.4 Known cases}
The results of \cite{Artin74} imply (1.1)  if the group
action is free. For varieties over $\c$ the conjecture was proved by
\cite{Popp73}, but his proof sheds no light on the algebraic nature of the
problem. He used  earlier results of \cite{Holman63; Kaup68} asserting that a
quotient exists as an analytic space, and then he  put an
algebraic space structure on it.

For seminormal complex analytic spaces Grauert developed a theory of
quotients by rather general equivalence relations; see
\cite{Dethloff-Grauert94} for a survey.  It would be interesting to work out
the analogous theory for algebraic spaces.

A simpler approach to the quotient problem over $\c$ was outlined in
\cite{Mumford-Fogarty82, p.172}. In its original form it still used complex
analytic quotients, but recently is has been reformulated in terms of
algebraic geometry \cite{Keel94; Viehweg95, 9.16}.  This method works  when
the stabilizers of points are reduced. The latter   always holds over $\c$, but
gives a  rather unnatural restriction  in general.
\enddemo

The aim of this article is to prove the existence of quotients in two cases.

\proclaim{1.5 Theorem} Fix an excellent base scheme $S$. Let  $G$ be an affine
algebraic group scheme   of finite type  over $S$ and $X$ a separated
algebraic space   of finite type  over $S$. Let $m:G\times X\to X$ be a proper
$G$-action on $X$. Assume that one of the following conditions is satisfied:

(1.5.1) $G$ is a reductive group scheme over $S$ (7.1).

(1.5.2) $S$ is the spectrum of a field of positive characteristic.
\smallpagebreak

\noindent Then a geometric quotient $f:X\to X/G$  exists, and $X/G$ is a
 separated
algebraic
space  of finite type over $S$.
\endproclaim

\demo{1.6 Comments}  (1.6.1) There are many examples of moduli problems where
the quotient of a scheme by a group is not a scheme. Thus it is natural to
consider the quotient problem for algebraic spaces.

(1.6.2) In most of the applications concerning moduli
problems, one has to take quotients by $PGL(n,S)$.  For these problems,
(1.5) provides a  satisfactory answer.

(1.6.3) It is somewhat surprising, but I do not see how to use my methods
when $S$ is the spectrum of a field of   characteristic zero and the group is
not reductive. Fortunately, this case has been settled earlier.
\enddemo

The general methods for constructing moduli spaces as geometric quotients
 imply that  (1.5) immediately yields the
existence of a slew of moduli spaces over $\spec \z$
(cf. \cite{Mumford-Fogarty82; Koll\'ar85}). These results were previously
known only in characteristic zero.

The results are most complete for the  arithmetic moduli of surfaces of
general type. Let
$\Cal M_{K^2,\chi}$ denote the moduli functor of canonical models  of
surfaces of general type with fixed
$K^2$ and
$\chi(\o)$.

\proclaim{1.7 Theorem} There is an algebraic space
$M_{K^2,\chi}$
over $\spec \z$ which coarsely represents the moduli functor
$\Cal M_{K^2,\chi}$.  $M_{K^2,\chi}$ is  separated and of finite type
over $\spec \z$. \qed
\endproclaim

For higher dimensional polarized varieties the situation is not  so fully
understood. Matsusaka's big theorem \cite{Matsusaka72} is still not known in
positive characteristic. There are further problems  with the locus of ruled
varieties (cf. \cite{Koll\'ar85, p.392}).

One   case where the latter problem does not occur is the following. Let
$\Cal M_{P}$ denote the moduli functor of smooth, canonically  polarized
varieties with Hilbert polynomial $P$.

\proclaim{1.8 Theorem} There is an algebraic space
$M_{P}$
over $\spec \z$ which coarsely represents the moduli functor
$\Cal M_{P}$.  $M_{P}$ is  separated and locally of finite type
over $\spec \z$.\qed
\endproclaim

\demo{1.9 Contents of the Sections}
Section 2 reworks the basic foundational properties of geometric quotients
in the category of algebraic spaces. The main reason necessitating this is
that among algebraic spaces a geometric quotient need not be unique (2.18).
Uniqueness holds for proper actions, but it needs a careful proof.

In section 3 I study general properties of quotient morphisms. The main result
(3.12) asserts that if a morphism looks like  a quotient topologically, then a
geometric quotient exists with the same underlying topological space.
These results are   related to \cite{Fogarty83}.

Next I study the quotient problem for successively larger classes of
algebraic spaces.

For normal spaces the quotient
problem was considered in \cite{Seshadri63,71}. He did not formulate the
results in terms of algebraic spaces, but in essence he proved (1.1) for $X$
normal. These theorems  are recalled in section 4.

If $X$ is reduced, we would like to construct the quotient using the
normalization $\bar X\to X$. We already know that $\bar X/G$ exists, and we
just have to undo the effects of normalization. These types of problems have
been considered by \cite{Artin70,6.3}.  The problem turns out to be rather
delicate. For technical reasons we
assume that $X$ is weakly normal (5.4).
The group theoretic aspects are discussed in section 5; most of the
technical problems are left to section 8.

For an arbitrary algebraic space the quotient is built  using
the already constructed quotient of the weak normalization of $X$. In positive
characteristic the Frobenius morphism can be exploited to give the general
result. This is discussed in section 6.

Reductive group schemes  are considered in section 7. I have to use the
theorems of
\cite{Seshadri77} on reductive group schemes acting on affine schemes, though
this goes against the general philosophy  of my approach.
\enddemo

\demo{1.10 Technical Assumptions} Throughout the paper, the base scheme is
denoted by
$S$. We always assume that $S$ is locally Noetherian. Starting with section 4,
we also assume that it is excellent. This assumption could be weakened a
little, but there are serious difficulties arising for general Noetherian base.

Starting with section 3, all algebraic spaces  are assumed separated and
locally of finite type over $S$.
Lack of separation  causes serious problems, but
they probably can be overcome through  careful examination of the proofs and
definitions. In considering group actions, the interesting part always seems
to happen in a  finite type extension, though I do not know any general result
to this effect. The assumpton is also important technically.

All group schemes are assumed to be separated and of finite type over $S$.
\enddemo

\demo{Acknowledgements}
 Partial financial support was provided by the
NSF under grant numbers DMS-8707320 and DMS-9102866.

These notes were typeset by \AmSTeX, the \TeX macro system of the American
Mathematical Society.
\enddemo

\head 2. Geometric Quotients of Algebraic Spaces
\endhead

The aim of this section is to review the basic definitions and theorems about
group actions and various quotients. Most of the material follows
\cite{Mumford-Fogarty82, chapter 0} rather closely,  reinterpreted in the
context of algebraic spaces.  Most definitons and results go through without
any change. The
main reason of repeating  the material for
algebraic spaces is that there is one subtle and surprising difference:

Geomeric quotients  are not necessarily categorical quotients any more, and
there may exist more than one geometric quotient (2.18).

For proper group actions these anomalies do not occur, but this fact requires
a careful proof (2.15).

\demo{2.1 Definition}  (2.1.1)  A {\it group scheme} over $S$ is
 a scheme $G/S$ together with morphisms

multiplication: $\mu:G\times_S G\to G$, inverse: $b:G\to G$,  and
unit: $e:S\to G$,

\noindent which satisfy the obvious group laws in the scheme theoretic setting
(see \cite{Mumford-Fogarty82, 0.1}).

(2.1.2) Let $G/S$ be a group scheme and $X/S$ an algebraic space.
An {\it action of $G$ on $X$} is a morphism
$m_X:G\times_SX\to X$ such that
$m(\mu(p_1,p_2),p_3)=m_X(p_1,m_X(p_2,p_3))$ where the $p_i$ are the coordinate
projections of $G\times_SG\times_SX$.
We usually say that $X$ is a {\it $G$-space} if the action is
understood.

If $\psi:U\to G$ and $f:U\to X$ are morphisms then $m_X(\psi, f):U\to X$ is
usually denoted by $\psi\cdot f$.

The {\it trivial} $G$ action on $X$ is given by $m_X=p_2$ where
$p_2:G\times_SX\to X$ is the second projection.

(2.1.3) Let $G/S$ be a group scheme, $X/S$ a $G$-space and  $F$ a sheaf on
$X$. An {\it action of
$G$ on $F$} is an isomorphism $\sigma:m_X^*F\cong p_2^*F$ which satisfies the
scheme theoretic analog of the cocycle condition; see
\cite{Mumford-Fogarty82, 1.3} for details.
\enddemo

\demo{2.1.4 Remark} One can  also  consider group objects in the category of
algebraic spaces, but these never seem to have acquired a good name. Aside
from this  terminological problem, we could always use this more
general notion.
\enddemo

\demo{2.2 Definition}  (2.2.1) Let $m_i:G\times X_i\to X_i$ be
$G$-spaces. A
morphism
$f:X_1\to X_2$ is called  a {\it $G$-morphism} if
the following  diagram is commutative:
$$
\CD
G\times X_1@>m_1>> X_1\\
@V{id_G\times f}VV @VV{f}V\\
G\times X_2@>m_2>> X_2.
\endCD
$$

If $F$ is a $G$-sheaf on $X_1$ then the  direct image $f_*F$ is
$G$-sheaf on $X_2$.

(2.2.2) A $G$-map of $G$-sheaves is defined in the obvious way. If $X$ is
an algebraic space with trivial $G$-action and $F$  a  (quasi coherent)
$G$-sheaf on $X$ then
$$
X\supset U\mapsto H^0(U,F|U)^G
$$
is a (quasi coherent) sheaf on $X$. It is denoted by $F^G$.
(The definiton makes sense even if the $G$-action on $X$ is not trivial, but
then $F^G$ is not a sheaf of $\o_X$-modules.)
\enddemo

\demo{2.3 Definition} Let $X$ be a $G$-space.  Let $\Psi_X$
denote the morphism
$$
\Psi_X:G\times_S X@>{(m_X,p_2)}>> X\times_S X,
$$
where $p_2:G\times_S X\to X$ is the second projection.

 We say that the  action of
$G$ on $X$ is  {\it proper} if  $\Psi_X$ is proper.
\enddemo

\demo{2.3.3 Remark}
It is somewhat unfortunate that the name ``proper" is
  used  to describe a condition which guarantees the
separatedness of the quotient (2.9).
\enddemo

The following is a reformulation of the valuative criterion of properness:

\proclaim{2.4 Lemma}   Let $X$ be a $G$-space.  The $G$ action on $X$ is
proper iff the following condition  is satisfied:

Let
$T$ be the spectrum of a DVR with generic point
$t_g$ and
$u_1,u_2:T\to X$ two morphisms. Assume that  there is a morphism $\psi_g:t_g\to
G$ such that $u_1(t_g)=\psi_g(t_g)\cdot  u_2(t_g)$.  Then
there is unique extension of  $\psi_g$ to a  morphism $\psi:T\to
G$ such that $u_1=\psi\cdot  u_2$.\qed
\endproclaim

\proclaim{2.5 Proposition}  Let $f:X\to Y$ be a  separated $G$-morphism of
$G$-spaces.
 If  $G$ acts properly on $Y$ then $G$ acts properly on $X$.
\endproclaim

\demop  By assumption $\Psi_Y: G\times Y@>>> Y\times Y$ is proper, By base
change we obtain that  $\Psi_Y\times_YX: G\times X@>>> Y\times X$ is also
proper. The latter can be factored as
$$
\Psi_Y\times_YX: G\times X@>\Psi_X>> X\times X@>{(f,id_X)}>> Y\times X.
$$
Thus $\Psi_X$ is proper.
\qed\enddemo

\demo{2.6 Definition}
Let $f_i:X_i\to Y$  be $G$-morphisms of $G$-spaces.
Let $\Delta:G\to G\times G$ be the diagonal imbedding. Then
the fiber product $X_1\times_YX_2$ is a
$G$-space via the morphism
$$
m:G\times(X_1\times X_2)@>{(\Delta,id)}>>G\times G\times(X_1\times X_2)
\cong (G\times X_1)\times(G\times X_2)@>{(m_1,m_2)}>> X_1\times X_2.
$$
The coordinate projections $X_1\times_YX_2\to X_i$
are $G$-morphisms.
\enddemo

There are various ways of defining the quotient of an algebraic space by a
group.  The following are the basic variants. I am mostly interested in
proper actions, therefore the concept of ``good quotient" is not relevant.

\demo{2.7 Definition} Let $G/S$ be a group scheme and $X/S$ a $G$-space.
Let $Z$ be an algebraic  space with trivial $G$-action and $q:X\to Z$ a
$G$-morphism.

(2.7.1) $q:X\to Z$ is called a {\it categorical quotient} of $X$ by $G$
if, for every algebraic space $U$ with trivial $G$-action and $G$-morphism
$p:X\to U$,  there is a morphism $u:Z\to U$  such that  $p=u\circ q$.

(2.7.2) $q:X\to Z$ is called a {\it topological quotient} of $X$ by $G$
if  the following   conditions are
satisfied:

(2.7.2.1) $f$ is a $G$-morphism.

(2.7.2.2) $f$ is locally of finite type.

(2.7.2.3) If $K$ is  an algebraically closed field then
$f(K):X(K)/G(K)\to Z(K)$ is an isomorphism (of sets).

(2.7.2.4) $f$ is universally submersive.

(2.7.3) $q:X\to Z$ is called a {\it geometric quotient} of $X$ by $G$
if  it is a topological quotient and in addition

(2.7.2.5) $\o_Z=(q_*\o_X)^G$.
\enddemo

\demo{2.8 Remarks}  (2.8.1)  Strictly speaking, the definitions (2.2--3)
make sense for algebraic spaces only if we establish their invarance under
\'etale morphisms $Z'\to Z$.  This is done in (2.10).

(2.8.2) By \cite{Mumford-Fogarty82, 0.1}, in the
category of {\it schemes} a geometric quotient is also a categorical quotient.
In  particular, it is unique up to isomorphism. We see in (2.18) that for
algebraic spaces both of these can fail if the $G$-action is not proper.

(2.8.3) This definition differs  slightly from the one in
\cite{Mumford-Fogarty82, sec. 0.1}. I assume that $f$ is locally of finite
type. In keeping with the original definition in \cite{Mumford65, sec. 0.1}, I
assume that $f$ is universally submersive (and not just submersive).
As \cite{Mumford65, sec. 0.1} argues, universally submersive is a more
topological notion. If $G$ is open over $S$ then the two notions are
equivalent by \cite{Dixmier-Raynaud81, 1.9}, and it is convenient to have the
stronger variant from the beginning.
For topological quotients it is crucial to assume universal submersiveness,
see (3.12.4).

(2.8.4) By \cite{Mumford-Fogarty82, p.6}, if $q$ is a topological quotient
and $G/S$ is universally open
then
$q$ is also universally open.
\enddemo

\proclaim{2.9 Proposition} Assume  that the $G$-action on $X$ is proper
and a topological quotient $f:X\to Z$ exists.
 Then $Z$ is separated.
\endproclaim

\demop    $f$ is universally submersive, thus, if $T$ is the
spectrum of a DVR and
$u:T\to Z$  a morphism then there is  a  dominant morphism $h:T'\to T$
(where
$T'$ is the spectrum of a DVR) such that $u\circ h:T'\to Z$ can be lifted to
a morphism $\bar u:T'\to X$ (cf. (3.7)).

Let $u_1,u_2:T\to X/G$ be two morphisms which agree at the generic point.
We find $T'\to T$ such that both of them lift to $\bar u_i:T'\to X$.
$\bar u_1|t'_g$ and $\bar u_2|t'_g$ lift the same morphism $t'_g\to Z$, thus
they differ by translation by $G$.  Since
the $G$ action on $X$ is proper, the two morphisms $\bar u_i$ differ only by
a translation by $G$. Thus $f\circ \bar u_1=f\circ \bar u_2$, hence $Z$ is
separated.
\qed\enddemo

\proclaim{2.10 Proposition} \cite{Mumford-Fogarty82, p.9} Let $X$ be a
$G$-scheme and $f:X\to Z$ a $G$-morphism (with the  trivial $G$-action on
$Z$). Let $h:Z'\to Z$ be a flat morphism. Set $X'\deq X\times_ZZ'$ and let
$f':X'\to Z'$  be the second projection.

(2.10.1) If $f:X\to Z$ is a
  geometric (resp. topological) quotient of  $X$ mod $G$, then
$f':X'\to Z'$ is also a geometric  (resp. topological) quotient of  $X'$ mod
$G$.

(2.10.2) Assume that $h$ is faithfully flat. Then  $f':X'\to Z'$ is a
  geometric (resp. topological) quotient of  $X'$ mod $G$ iff
$f:X\to Z$ is   a geometric  (resp. topological) quotient of  $X$ mod $G$.\qed
\endproclaim

If $f:X\to Y$ is an \'etale $G$-morphism then the induced morphism
between the geometric quotients is usually not \'etale. The following
condition was introduced by Deligne (cf. \cite{Knutson71, p. 183})
to remedy this problem.

\demo{2.11 Definition}  Let $X$ be a $G$-space, and $\Psi_X$ as in (2.3).
The fiber   $\Psi_X^{-1}(x,x)\subset G\times \{x\}$ is a subgroup scheme. It
is  called the {\it stabilizer} of $x$, and is denoted by $\stab(x)$.
\enddemo

\demo{2.12 Definition}  Let $f:X\to Y$ be a
$G$-morphism. We say that $f$ is {\it fixed-point reflecting}  if the
following condition is satisfied:

$\stab(x)=\stab(y)\times_{y}\{x\}$ for  any
$x\in X$  and $y=f(x)$.
\enddemo

The simplest examples of fixed-point reflecting morphisms are obtained by
base change from trivial $G$-actions:

\proclaim{2.13 Lemma} Let $f:X\to Z$ be a $G$-morphism with trivial $G$-action
on $Z$ and $h:Z'\to Z$ any morphism. Then $X\times_ZZ'\to X$ is
fixed-point reflecting.\qed
\endproclaim

The following result shows the importance of
fixed-point reflecting morphisms.
It is related to \cite{Luna73, Chap. II}.
The theorem is also crucial in establishing
basic properties of geometric quotients in the category of algebraic spaces.

\proclaim{2.14 Theorem} Assume that $G\to S$ is universally open.
Let $f:X\to Y$ be a separated and fixed-point
reflecting
$G$-morphism.  Assume that the $G$-action on $Y$ is proper and that
$f$ is
\'etale. Assume furthermore that a topological quotient $p_Y:Y\to Z(Y)$ exists.

Then there is a topological quotient $p_X:X\to Z(X)$ such that the
induced morphism
$Z(X)\to Z(Y)$ is
\'etale.
\endproclaim

\demop  By (2.8.4)  $p_Y$ is universally open and surjective. Thus by
\cite{SGA1, IX.4.9}  all descent data on separated and \'etale morphisms
to
$Y$ are effective. (See \cite{BLR90, chapter 6} for a clear and selfcontained
discussion on descent.) Therefore it is sufficient to construct a descent
datum for $f:X\to Y$, that is, an isomorphism
$$
\phi:X\times_{Z(Y)}Y@>\cong>> Y\times_{Z(Y)}X\qtq{over} Y\times_{Z(Y)}Y
$$
which satisfies the cocycle condition \cite{BLR90, p.133}.

The construction of $\phi$, and the role of the fixed-point reflecting
conditions, is   transparent at the set theoretic level.

Let $K$ be an algebraically closed field.
Pick $x\in X(K)$ and $y\in Y(K)$ such that $p_Y(f(x))=p_Y(y)$. By (2.7.2.3)
there is a
$g\in G(K)$ such that $f(x)=g\cdot y$. Set
$$
\phi(x,y)=(g\cdot y, g^{-1}\cdot x).\tag 2.14.1
$$
$g$ is not uniquely defined by the equation $f(x)=g\cdot y$; two choices
differ by an element of $\stab(y)$. Since $\stab(x)=\stab(y)$, the product
$g^{-1}\cdot x$ is independent of the choice of $g$. Thus $\phi$ is well
defined.

We construct $\phi$ by constructing its graph.  (2.14.1) shows
 that, at least set theoretically,  the graph of $\phi$ is the
image of the morphism
$$
\align
&h:G\times_{Z(Y)}X\to
V\deq X\times_{Z(Y)}Y\times_{Y\times_{Z(Y)}Y}Y\times_{Z(Y)}X,
\qtq{given by}\\
&h:(g,x)\mapsto (g\cdot x, f(x), g\cdot f(x),  x ).
\endalign
$$
$G\times_{Z(Y)}Y\to Y\times_{Z(Y)}Y$ is proper, hence so is
$h_2: G\times_{Z(Y)}X\to Y\times_{Z(Y)}X$. $h$ is the fiber product of two such
morphisms, hence $h$ is also proper. Thus $W\deq \im h\subset V$ is a closed
subscheme.

$h_2$ can be factored as
$$
h_2:G\times_{Z(Y)}X@>>> W @>\pi_2>> Y\times_{Z(Y)}X.
$$
This shows that $\pi_2$ is proper. The set theoretic
computation of $\phi$ shows that $\pi_2$ is universally injective and
surjective, thus $\pi_2$ is a universal homeomorphism.

$X\times_{Z(Y)}Y$ and $Y\times_{Z(Y)}X$ are both \'etale over
$Y\times_{Z(Y)}Y$, thus  the projection morphisms
$$
\CD
V@>>> Y\times_{Z(Y)}X\\
@VVV @VV{id_Y\times f}V\\
X\times_{Z(Y)}Y@>{f\times id_Y}>> Y\times_{Z(Y)}Y
\endCD
$$
are all \'etale.

Thus $W\subset V$ provides a section of $V@>>> Y\times_{Z(Y)}Y$, at least
after  base change  to $W$. Since
$W @>>> Y\times_{Z(Y)}X$ is  a universal homeomorphism,   (5.3) implies that
$V@>>> Y\times_{Z(Y)}Y$ has a unique section $W^0\subset V$ such that $\supp
W^0=\supp W$. $W^0$ is the  graph of  $\phi$.

By (5.3) it is sufficient to check the cocycle condition  set
theoretically, where it is trivial.\qed\enddemo

\demo{2.14.2 Remark} (2.14) implies that a fixed-point reflecting and \'etale
$G$-morphism is also strongly \'etale in the sense of \cite{Mumford-Fogarty82,
p. 152}.
\enddemo

\proclaim{2.15 Corollary} Let $X$ be an algebraic $G$-space. Assume that the
$G$-action on $X$ is proper. Let $q:X\to Z$ be a geometric quotient.
Then $q$ is also a categorical quotient. Thus a geometric quotient is unique
(up to a canonical isomorphism).
\endproclaim

\demop
\cite{Mumford-Fogarty82, p.4} shows that in the
category of {\it schemes} a geometric quotient is also a categorical quotient.
 More generally, the proof
shows the following:

\proclaim{2.15.1 Claim} Let $X$ be an algebraic $G$-space and $q:X\to Z$  a
geometric quotient.  Let $V$ be a scheme with trivial $G$-action and
$p:X\to V$  a $G$-morphism. Then there is a morphism $v:Z\to V$ such that
$p=v\circ q$.\qed
\endproclaim

In order to get the general case, let $U$ be an algebraic space with trivial
$G$-action and
$p:X\to U$  a $G$-morphism. We need to find $u:Z\to U$. To do this let $V$ be
an affine scheme and $j:V\to U$ an \'etale morphism. By base change we obtain
a fixed-point reflecting
\'etale morphism $p^*j:V\times_UX\to X$.

Applying (2.14) we obtain that a geometric quotient
$p':V\times_UX\to Z'$ exists and the induced morphism $Z'\to Z$ is \'etale.
Since $V$ is a scheme, by (2.15.1) we obtain a morphism $Z'\to V$.
The resulting diagram of \'etale morphisms
$$
Z@<<< Z'@>>> V@>>> U
$$
gives the local description of  $u:Z\to U$.
\qed\enddemo

\demo{2.16 Notation} Let $X$ be an algebraic $G$-space such  that the
$G$-action on $X$ is proper. The unique geometric quotient is denoted by
$q:X\to X/G$.
\enddemo

\proclaim{2.17 Corollary} Let $f:X\to Y$ be a separated, surjective, \'etale
and fixed-point reflecting
$G$-morphism.  Assume that the $G$-action on $Y$ is proper. Assume furthermore
that the geometric  quotient $q_X:X\to X/G$ exists.

Then the geometric quotient $q_Y:Y\to Y/G$ exists and the
induced morphism
$X/G\to Y/G$ is
\'etale.
\endproclaim

\demop Let $Z=X\times_YX$. $Z$ can be viewed as an \'etale equivalence
relation on $X$ and  the quotient of $X$ by this equivalence
relation is $Y$.

Apply (2.14) to $Z\to X$. We conclude that $Z/G$ exists and the maps
$Z/G\to X/G$ are \'etale. Thus $Z/G$ gives an \'etale equivalence
relation on $X/G$. The quotient is an algebraic space $U$ such that
$Y\times_U(X/G)\cong X$. Thus by (2.10.2) $U$ is the
geometric quotient of $Y$ mod $G$.\qed\enddemo

The following example, taken form \cite{Mumford-Fogarty82, p.11}
shows that  a geometric quotient need not be a categorical quotient in general.

\demo{2.18 Example} Let $SL(2,\c)$ act on homogeneous polynomials by the rule
$$
\left(
\matrix
a&b\\
c&d\\
\endmatrix
\right): x\mapsto ax+cy,\ y\mapsto bx+dy.
$$
Let $X$ be the variety of pairs $(L,Q^2)$
and $Y$ the variety of pairs $(L,Q)$
where $L$ is a nonzero linear form
and $Q$ is a quadratic form of discriminant -1. It is easy to see that $X$
and $Y$ are smooth, irreducible  and $SL(2,\c)$ acts freely on them. The
morphism
$$
f:Y\to X\qtq{given by} (L,Q)\mapsto (L,Q^2)
$$
is an \'etale and fixed-point reflecting $SL(2,\c)$-morphism. $f$ has degree
two and the fiber over $(L,Q^2)$ has two points $(L,Q)$ and $(L,-Q)$.

A geometric quotient of $X$ is obtained as follows. Every pair can be
brought to the form $(x,Q^2)$. Then
$$
X\to \c\qtq{given by}(x,Q^2)\mapsto Q(0,1)^2
$$
is a geometric quotient.

Since $\c$ has no nontrivial \'etale covers, this already shows that (2.14)
can fail if the group action is not proper.

The stabilizer of $x$ is the group of upper triangular matrices
$$
\left(
\matrix
1&b\\
0&1\\
\endmatrix
\right).
$$
On quadratic forms the action of upper triangular matrices is
$$
ux^2+vxy+wy^2\mapsto (u+vb+wb^2)x^2+(v+2wb)xy+wy^2.
$$
This shows that $(L,Q)$ and $(L,-Q)$ are always in different orbits.

A geometric quotient of $Y$ is obtained as follows. Every pair can be
brought to the form $(x,Q)$. Then
$$
X\to Z\cong \c\qtq{given by}(x,Q)\mapsto Q(0,1)
$$
is a geometric quotient, provided $Q(0,1)\neq 0$. If $Q(0,1)=0$ then we get
two different orbits corresponding to the pairs $(x,xy)$ and $(x,-xy)$.
The geometric quotient is not separated, as shown by the families of pairs
$$
(x,xy+ty^2)\qtq{and} (x,-xy+ty^2).
$$
These lie in the same orbit for $t\neq 0$ but in different orbits if $t=0$.
Thus a geometric quotient  of $Y$ by $SL(2,\c)$ is the nonseparated line
(which is a scheme)
$$
U:\ \c\times\{0,1\}\qtq{modulo the equivalence} (t,0)\sim_1 (t,1)\qtq{for
$t\neq 0$.}
$$
The \'etale equivalence relation $(t,0)\sim_2 (-t,1)$ on
$\c\times\{0,1\}$ descends to an \'etale equivalence relation on $U$.
Let $V$ denote the quotient (which is an algebraic space).

The  morphism
$$
\c\times\{0,1\}\to  Z\cong \c\qtq{given by}
(t,0)\mapsto t^2,\ (t,1)\mapsto t^2
$$
descends to a morphsim $p:V\to Z$. $p$ is an isomorphism on the set of points,
but it is not an isomorphism. Indeed,  $\partial/\partial t$
gives a nonzero tangent vector of $V$ at the origin, but it maps to zero on
$Z$.
$V$ is  the simplest example of a locally nonseparated algebraic
space (cf. \cite{Knutson71, p.9}).

Finally, $Y\to U$ descends to an $SL(2,\c)$-morphism $X\to V$, but
there is no corresponding morphism from $Z$ to $V$. Thus the geometric
quotient is not a categorical quotient in the category of algebraic spaces.

One can see that in this case $V$ is both a geometric and
 a categorical quotient of $X$.
\enddemo

\demo{2.19 Example} Let $S=\c$ and  $G_1\deq S\times \{+1,-1\}$ the
constant 2-element group scheme over $S$. Let $X\subset \c\times S$ be the
union of the diagonal $(x,x)$ and of the anti diagonal $(-x,x)$. $G_1$ acts
on $X$ by multiplication on the first coordinate.

Let $G\subset G_1$ be the subgroup scheme $S\times \{+1\}\cup
(S\setminus\{0\}) \times
\{-1\}$. $G$ acts on $X$, the action is closed but not proper and the
geometric quotient is $S$. $X\to S$ is affine; this is of interest in
connection with (3.12).

The $G$-action on $X$ does not identify the two directions of the tangent
cone of $X$. One can see that the categorical quotient (in the category of
algebraic spaces) is not locally separated above $0\in S$.
\enddemo

\demo{2.20 Remark} It is quite likely that a geometric quotient is also a
categorical  quotient  in the category  of locally separated algebraic spaces.
If  $S=\c$ and $G$ is reductive,   this follows from the fact that a geometric
quotient is also a categorical  quotient  in the category of complex analytic
spaces
\cite{Neeman89, Thm. 10}.
\enddemo

\head 3. Topological Quotients
\endhead

Let $X$ be a $G$-space.   In this
section  we assume that we already have a quotient morphism $f:X\to Z$, and
concentrate on analyzing the relationships among its properties.   The main
result is (3.12), which asserts that a quotient with good topological
properties has very good scheme theoretic properties as well.

\demo{3.1 Assumptions} Fix a locally Noetherian base scheme $S$. For the rest
of this section we assume that
$G$ is  locally of finite type over $S$ and the structure morphism $G\to S$ is
universally open.

All algebraic spaces are assumed to be separated. All group
actions are assumed to be proper.
\enddemo

\demo{3.2 Definition} Let $X$ be a $G$-space. A morphism $f:X\to Z$ is
called a {\it weak  topological quotient} if the following conditions are
satisfied:

(3.2.1) $f$ is a $G$-morphism (with the trivial $G$-action on $Z$).

(3.2.2) $f$ is locally of finite type.

(3.2.3) If $K$ is  an algebraically closed field then
$f(K):X(K)/G(K)\to Z(K)$ is an isomorphism (of sets).
\enddemo

If $Z$ is normal then these properties are quite strong (3.9).
The
following example shows  what can go wrong in general:

\demo{3.3 Examples}  Let $S=\spec \c$ and $G=S$  the trivial group scheme.

(3.3.1)  Let $Z=\spec \c[x]$,  $X=\{0\}\cup (Z\setminus \{0\})$ (disjoint
union) and $f:X\to Z$ the natural morphism. The conditions (3.2.1--3)
are satisfied.  A more sophisticated version of this is example is the
following.

(3.3.2)  Let
$Z$ be an irreducible  curve with a single node $z\in Z$ and  $f':X'\to Z$ the
normalization with
$x_1,x_2\in X'$ the preimages of $z$. Set $X\deq X'-x_2$ and $f:X\to Z$ the
restriction of $f'$. Then $f$ is a homeomorphism and the conditions (3.2.1--3)
are satisfied.

Two bad properties of $f$ come to mind first: $f$ is not proper and not flat.
Unfortunately,     in  most cases, the quotient morphism is neither flat
nor proper.

The correct approach seems to find suitable variants of the valuative
criterion of properness. The next two definitions give first a  weak
variant and then a much stronger one. At the end, they turn out to be
equivalent (3.7).
\enddemo

\demo{ 3.4 Definition}  Let $f:X\to Y$ be a morphism of algebraic spaces.
We say that $f$ satisfies the
{\it   weak lifting property for DVR's}  if  the following holds:

Let $T$ be the spectrum of a DVR and $u:T\to Y$ a morphism.
 Then there is a commutative diagram
$$
\CD
T'@>u'>> X\\
@VhVV @VVfV\\
T@>u>> Y,
\endCD
$$
where $T'$ is the spectrum of a DVR  and
$h:T'\to T$ a  surjective  morphism.
\enddemo

\demo{ 3.5 Definition}  (cf. \cite{Seshadri72, 4.1})  Let $f:X\to Y$ be a
$G$-morphism of
$G$-spaces such that the $G$-action on $Y$ is trivial.

We say that $f$ is
{\it   complete  mod $G$}  if $f$ is separated,
of finite type and the following version of the valuative criterion of
properness is satisfied:

Let $T$ be the spectrum of a DVR with generic point $j:t_g\to T$.  Assume
that we have a commutative diagram
$$
\CD
t_g@>u_g>> X\\
@VjVV @VVfV\\
T@>>> Y
\endCD\tag 3.5.1
$$
 Then there is a commutative diagram
$$
\CD
T'@>u'>> X\\
@VhVV @VVfV\\
T@>>> Y
\endCD
\quad\qtq{and}\
\matrix
\text{a morphism}\\
\phi'_g:t'_g\to G,
\endmatrix
\tag 3.5.2
$$
where $T'$ is the spectrum of a DVR with generic point $t'_g$ and
$h:T'\to T$ a  dominant  morphism,
 such that
$$
u_g\circ (h|t'_g)=\phi'_g\cdot (u'|t'_g).
$$

Informally speaking, the valuative criterion of properness is satisfied
after a base change and a suitable translation by the group.
\enddemo

\demo{3.5.3 Remarks}  (3.5.3.1) If the $G$-action on $X$ is trivial
then $f$  is proper  iff it is complete mod $G$.

(3.5.3.2) Since $f$ is locally of finite type, we may assume in addition that
$h:T'\to T$ is finite.

(3.5.3.3) It is possible to define the notion of ``complete mod $G$" for
$G$-morphisms where we do not assume that the $G$-action on $Y$ is trivial. I
do not see any advantages at the moment.
\enddemo

The following proposition  lists some easy results about the composition of
morphisms which are proper or complete mod $G$.  I list only those  variants
which are used later. The proofs are straightforward.

\proclaim{3.6 Proposition}  Let $f:X\to Y$ and $g:Y\to Z$ be
separated $G$-morphisms of
$G$-spaces.

(3.6.1) If $f$ is proper and $g$ is complete mod $G$ then $g\circ f$ is
complete mod $G$.

(3.6.2) If $f$ is complete mod $G$ and $g$ is  proper then $g\circ f$ is
complete mod $G$.

(3.6.3) If $g\circ f$ is
complete mod $G$  and the $G$ action on $Y$ is trivial then
$f$ is complete mod $G$.

(3.6.4) If $g\circ f$ is
complete mod $G$, the $G$ action on $Y$ is trivial and $f$ is surjective then
$g$ is proper.\qed
\endproclaim

The next result shows that several   good properties of
weak topological quotients are equivalent:

\proclaim{3.7 Proposition}   Let $X$ be a $G$-space and  $f:X\to Z$  a
$G$-morphism, locally  of finite type (with trivial $G$-action on $Z$). Then
each of the  following conditions implies the next:

(3.7.1) $f$  is surjective and
complete mod $G$.

(3.7.2) $f$ is surjective  and universally open.

(3.7.3) $f$ is universally submersive.

(3.7.4) $f$ satisfies the   weak lifting property for DVR's.
\smallpagebreak

If, in adition, $f:X\to Z$  is a
weak topological quotient, then all four condition are equivalent.
\endproclaim

\demop Assume (3.7.1). (3.7.1) is stable under any base change, thus it is
sufficient to prove that $f$ is open. Let $U\subset X$ be open and
$V=f(U)\subset Y$. $V$ is constructible by Chevalley's theorem \cite{EGA71,
I.7.1.4}, hence
$V$ is open iff it is closed under generalization.

Pick $y\in V$ and let $y'$ be a generalization of $y$. Let $h:T\to Y$ be a
morphism such that $y'$ is the
image of the generic point and $y$ the
image of the closed point.
$f$ is surjective, hence there are  $x,x'\in X$ such that
$f(x)=y$ and $f(x')=y'$.  By (3.5) there is a commutative diagram
$$
\CD
T'@>u'>> X\\
@VhVV @VVfV\\
T@>u>> Y.
\endCD
$$
By enlarging $T'$ we may assume that $T'$ has algebraically closed residue
field $K$ and that
$x$ is a $K$-valued point of $X$.  Thus there is a $\phi_0:K\to G$ such that
$x=\phi_0\cdot (u'|t'_0)$. After replacing $T'$ by a ramified extension,  we
may also assume that $\phi_0$ extends to $\phi:T'\to G$ (here we use that
$G/S$ is universally open). Set
$u''\deq
\phi\cdot u'$.

 $u'':T'\to X$ is a morphism such that  $x=u''(t'_0)$ and $x''\deq u''(t'_g)$
is a generalization of $x$. Thus $x''\in U$ and therefore $y'=f(x'')\in V$.
Thus $f$ is open.

(3.7.2) $\Rightarrow$ (3.7.3) holds for any continuous map of
topological spaces.

In order to see (3.7.3) $\Rightarrow$ (3.7.4)  let $u:T\to Z$ be a morphism
and
$t_0\in T$ the closed point. By assumption $u^*f:X\times_ZT\to T$ is
submersive, thus
$(u^*f)^{-1}(t_0)\subset X\times_ZT$ is not open. Thus there is a point
$x\in (u^*f)^{-1}(t_0)$ which is the specialization of a point $x'\not\in
(u^*f)^{-1}(t_0)$.  This gives  $v:T'\to X\times_ZT$  such that $x'$ is the
image of the generic point and $x$ the
image of the closed point. Thus $h:T'\to T$ is dominant and
$u':T'\to X\times_ZT\to X$ gives a lifting.

Finally,  assume (3.7.4) and that $f$ is a  weak topological
quotient. Consider a commutative diagram
$$
\CD
t_g@>u_g>> X\\
@VjVV @VVfV\\
T@>u>> Y.
\endCD
$$
By assumption there is a commutative diagram as in (3.4)
$$
\CD
T'@>u'>> X\\
@VhVV @VVfV\\
T@>u>> Y.
\endCD\tag 3.7.5
$$
By enlarging $T'$ we may assume that the quotient field $K$ of $T'$ is
algebraically closed. We show that (3.7.5) is the diagram whose existence is
required in (3.5.2).

$u'|t'_g$ and $u_g$ are $K$-valued points of $X$ and they have the same
image in $Y$. Since  $f$ is a  weak topological
quotient,  there is a $\phi'_g:t'_g\to G$ such that
$u_g\circ (h|t'_g)=\phi'_g\cdot (u'|t'_g)$.  This shows
 (3.7.1).
\qed\enddemo

\demo{3.7.6 Remark} The above proof shows that if $f:X\to Y$ is a morphism of
finite type then
 $f$ is universally submersive
iff $f$ satisfies the   weak lifting property for DVR's.
\enddemo

We can reformulate the definition (2.7.2)  as follows:

\demo{3.8 Definition}  Let $X$ be a $G$-space. $f:X\to Z$ is   a {\it
topological quotient} iff  $f$ is a weak topological quotient and $f$ satisfies
the equivalent conditions of (3.7).
\enddemo

If $Z$ is normal, this notion does not give anything new:

\proclaim{3.9 Proposition} Let $f:X\to Z$ be a  weak topological quotient.
Assume that
$Z$ is normal, irreducible  and every irreducible component of $X$ dominates
$Z$. Then
$f:X\to Z$ is a topological quotient.
\endproclaim

\demop  We want to use the Chevalley criterion   \cite {EGA, IV.14.4.4}  to
show that
$f$ is universally open. We need to prove that the fiber
dimension  $x\mapsto \dim_x f^{-1}(f(x))$ is locally constant on $X$.
Let $g:X\to S$ be the   morphism  giving the $S$-space structure. For
$s\in S$ let
$G_s$ denote the fiber of $G$ over $s$.

$x\mapsto \dim_x f^{-1}(f(x))$ and $x\mapsto \dim \stab(x)$  are  always upper
semi continuous. Furthermore,
$$
\dim_x f^{-1}(f(x))+\dim \stab(x)=\dim G_{g(x)}
$$
is locally constant since $G\to S$ is open. Thus $x\mapsto \dim_x f^{-1}(f(x))$
is locally constant on
$X$.\qed\enddemo

We introduce a further notion, whose main point is to allow quotients
 $f:X\to Z$ where
$f(K):X(K)/G(K)\to Z(K)$ is only finite--to--one. For such morphisms the
conditions in (3.7) are no longer equivalent. We retain the strongest one of
them.

\demo{3.10 Definition}  Let $X$ be a $G$-space. $f:X\to Z$ is called an {\it
approximate quotient} if the following conditions are satisfied:

(3.10.1) $f$ is a $G$-morphism
 (with the trivial $G$-action on $Z$).

(3.10.2) $f$ is of finite type.

(3.10.3) If $K$ is  an algebraically closed field then
$f(K):X(K)/G(K)\to Z(K)$ is surjective and finite--to--one (as a map of
sets).

(3.10.4) $f$ is complete mod $G$.
\enddemo

Approximate quotients are stable under many operations:

\proclaim{3.11 Proposition} Let $f:X\to Y$ and $g:Y\to Z$ be $G$-morphisms.

(3.11.1) Assume that $g$ is finite, surjective and the $G$-action on $Y$ is
trivial. Then $f$ is an approximate quotient iff $g\circ f$ is.

(3.11.2) Assume that $f$ is finite and surjective. Then $g$ is an approximate
quotient iff $g\circ f$ is.

(3.11.3)  Assume that  $f$ is an approximate (resp. topological) quotient and
let
$h:Y'\to Y$ be a morphism. Then $h\times f:Y'\times_YX\to Y'$ is an approximate
 (resp. topological) quotient.

(3.11.4) Let $\phi:H\to G$ be a homomorphism of algebraic groups which is
finite and maps onto a finite index subgroup of $G$. $\phi$ induces an
$H$-action on $X$. In this case, if $f$ is an approximate quotient by $G$ then
$f$ is an approximate quotient by $H$. The converse also holds if the
$G$-action on $Y$ is trivial.
\endproclaim

\demop The last assertion is clear.
In the remaining three cases the conditions (3.10.1--3) are easy to check.
In the first two cases (3.6) allows us to check (3.10.4), while this is
clear for (3.11.3). \qed\enddemo

The following is the main theorem of the section. It shows that approximate
quotients have surprisingly good properties.
The result can be viewed as a generalization of \cite{Fogarty83}.
In the complex analytic case similar questions are considered in
\cite{Roberts86}.

\proclaim{3.12 Theorem} Let $G$ be an affine algebraic group scheme,
universally open and locally of finite type over $S$. Let $m:G\times X\to X$
be a proper $G$-action on an algebraic space $X$.
 Let  $f:X\to Z$ be an approximate quotient.
 Then:

(3.12.1) $f$ is affine.

(3.12.2) If $X$ is of finite type over $S$ then so is $Z$.

(3.12.3)  Assume in addition that $G$ is flat over $S$. If $F$ is any coherent
$G$-sheaf on
$X$ then
$(f_*F)^G$ is a coherent sheaf on $Z$.
\endproclaim

\demo{3.12.4 Remark} \cite{Fogarty83} asserts a version of this assuming only
that  $f$ is submersive. Example (3.3) is a strict quotient in his sense
where his assertion about coherence fails. His assertion and proof remain
valid for $f$ universally submersive.
\enddemo

The following   consequence of (3.12)  shows that topological quotients
are very close to geometric quotients:

\proclaim{3.13 Theorem} Let $G$ be an affine algebraic group scheme,
flat and locally of finite type over $S$. Let $m:G\times X\to X$
be a proper $G$-action on an algebraic space $X$.
 Let  $f:X\to Z$ be a topological quotient.
Set
$$
X/G\deq \spec_Z(f_*\o_X)^G,\qtq{and let }f:X@>g>> X/G@>q>> Z
$$
be the natural morphisms.
Then $g:X\to  X/G$ is the geometric quotient of $X$ mod
$G$ and $q:X/G\to Z$ is a finite and universal homeomorphism (5.1).
\endproclaim

\demop By (3.12.3) $(f_*\o_X)^G$ is a coherent sheaf of algebras on $Z$,
thus $q:X/G\to Z$ is finite and surjective.

$g:X\to X/G$ is a $G$-morphism and the $G$-action on $X/G$ is trivial. $g$ is
locally of finite type since  $f$ is.
If $K$ is an algebraically closed field,   the composite
$$
X(K)/G(K)\to (X/G)(K)\to Z(K)
$$
is an isomorphism and $(X/G)(K)\to Z(K)$ is surjective. Thus
$X(K)/G(K)\to (X/G)(K)$ is an isomorphism.
By (3.11.1) $g$ is an approximate quotient, hence $g$ is universally
submersive.

Finally, $(g_*\o_X)^G=\o_{X/G}$ by constrction. Thus $g:X\to X/G$ is
the geometric quotient of $X$ mod
$G$.\qed\enddemo

\demo{3.13.1 Remark} In (3.13)  assume that $f$ is only an approximate
quotient. Then $g$ satisfies all the properties of geometric quotients except
possibly (2.7.2.3). I hope to study this problem later.
\enddemo

\demo{3.14 Proof of (3.12)}    First we show that
the properties (3.12.1--3) descend through finite and dominant morphisms.

The properties (3.12.2) and (3.12.1) can be studied for any morphism of
schemes.

\proclaim{3.14.1 Lemma} \cite{Atiyah-MacDonald69, 7.8; EGA71, 0.6.4.9} Let
$q:U\to V$ be a finite and surjective morphism of algebraic spaces. Then

(3.14.1.1) $U$
is (locally) of finite type over
$S$ iff $V$ is.

(3.14.1.2) $U$
is Noetherian  iff $V$ is.\qed
\endproclaim

The following is the relative version of the Chevalley theorem:

\proclaim{3.14.2 Lemma} \cite{EGA71, I.7.1.4; Knutson71, p.169} Consider a
commutative
 diagram  of algebraic spaces
$$
\CD
U'@>q_U>> U\\
@Vf'VV @VVfV\\
V'@>q_V>>V,
\endCD
$$
where $q_U$ and $q_V$ are finite and surjective. Then $f$ is affine iff $f'$
is.\qed
\endproclaim

The descent of coherence requires  a little care:

\proclaim{3.14.3 Lemma}   Let $G$ be an affine algebraic group
scheme, flat and of finite type over $S$. Let
 $m:G\times U\to U$ resp.  $m':G\times U'\to U'$ be  proper actions with
approximate quotients $f:U\to V$ resp. $f':U'\to V'$.
Assume that we have  a commutative
 diagram
$$
\CD
U'@>q_U>> U\\
@Vf'VV @VVfV\\
V'@>q_V>>V,
\endCD
$$
where $q_U$ and $q_V$ are finite and surjective and $q_U$ is a $G$-morphism.
Then the following are equivalent:

(3.14.3.1) If $F$ is any coherent $G$-sheaf on $U$ then $(f_*F)^G$ is a
coherent sheaf on $V$.

(3.14.3.2) If $F'$ is any coherent $G$-sheaf on $U'$ then $(f'_*F')^G$ is
a coherent sheaf on $V'$.
\endproclaim

\demop  If $F'$ is a coherent $G$-sheaf on $U'$ then $(q_U)_*F'$ is
a coherent $G$-sheaf on $U$ and
$$
(q_V)_*\left((f'_*F')^G\right)=\left(f_*((q_U)_*F')\right)^G.
$$
This shows that (3.14.3.1) implies (3.14.3.2).

The converse is proved by induction (with respect to lexicographic ordering)
on the pair
$(\dim\supp F, q(F))$ where $q(F)\deq \sum_i \len_{x_i}F$ where the $x_i$ are
the $\dim\supp F$-dimensional generic points of $\supp F$.

We need the following result:

\proclaim{3.14.3.3 Claim} \cite{SGA3, I.5.3; Seshadri77, p.264}
Assume that $G$ is flat over $S$.
Let $0\to
F_1\to F_2\to F_3$ be an exact sequence of quasicoherent $G$-sheaves on an
algebraic space with trivial
$G$-action.  Then the sequence
$$
0\to  F_1^G\to F_2^G@>>> F_3^G
$$
is exact.\qed
\endproclaim

Let $\frak S$ be the class of coherent $G$-sheaves $F$ such that $(f_*F)^G$
is coherent. We need to prove that $\frak S$ is the class of all coherent
$G$-sheaves.

\proclaim{3.14.3.4 Claim} Let $0\to F_1\to F_2\to F_3$ be an exact
sequence of coherent $G$-sheaves on $U$. If  $F_1,F_3\in \frak S$  then
$F_2\in \frak S$.
\endproclaim

\demop
By (3.14.3.3)  we have an exact sequence
$$
0\to (f_*F_1)^G\to (f_*F_2)^G@>v>> (f_*F_3)^G.
$$
Set  $E=\im v$. If $(f_*F_3)^G$ is coherent then $E\subset (f_*F_3)^G$
is also coherent. If in addition $(f_*F_1)^G$ is coherent then
$(f_*F_2)^G$ is also coherent.\qed
\enddemo

Let $F$ be a coherent $G$-sheaf on $U$. By induction,
if $E$ is a coherent $G$-sheaf on $U$ such that
$(\dim\supp E,q(E))<(\dim\supp F,q(F))$ then $E\in \frak S$.
We need to  prove that $F\in \frak S$.

Let $E$ denote the kernel of the natural map  $p:F\to (q_U)_*(q_U^*F)$.   We
have an exact sequence
$$
0\to E\to F@>p>> (q_U)_*(q_U^*F).
$$
$p$  is
nonzero at all generic points of
$\supp F$, thus $(\dim\supp E,q(E))<(\dim\supp F,q(F))$, hence $E\in \frak S$.
Set $F'\deq q_U^*F$. Then
$$
(q_V)_*\left((f'_*F')^G\right)=\left(f_*((q_U)_*F')\right)^G=
\left(f_*((q_U)_*(q_U^*F))\right)^G.
$$
By (3.14.3.1) $(f'_*F')^G$ is coherent on $V'$, thus
$(q_V)_*\left((f'_*F')^G\right)$ is coherent. This shows that
$(q_U)_*(q_U^*F)\in \frak S$.
By (3.14.3.4)  this implies that $(f_*F)^G$ is coherent.\qed\enddemo

The rest of the proof of (3.12)  relies on a slight strengthening of a result
of Mumford (see
\cite{EGA IV.14.5.10}):

\proclaim{3.14.4 Proposition}  Let $f:X\to Z$ be a universally open and
surjective morphism, locally of finite type. Then there is a finite and
surjective morphism $g:Z'\to Z$ and an open cover $\cup Z'_i=Z'$ such that
for every irreducible component $X'_{ij}\subset X'_i\deq X\times_ZZ'_i$
the induced morphism
$$
f'_{ij}:X'_{ij}\to Z'_i\qtq{has  a section} s_{ij}:  Z'_i\to X'_{ij} \qtq{for
every
$i,j$.}\qed
$$
\endproclaim

Apply (3.14.4) to $f:X\to Z$ in (3.12) and set $X'\deq X\times_ZZ'$.
$f':X'\to Z'$ is an aproximate quotient by (3.11.3). By (3.14.1--3) it is
sufficient to prove (3.12) for $f':X'\to Z'$.

The conditions (3.12.1--3) are local on $Z'$, thus we may replace $Z'$ by
$Z'_i$ to obtain $f'_i:X'_i\to Z'_i$. Let
$X'_{ij}\subset X'_i,\ j=1,\dots, j(i)$
be the irreducible components.  For every $i,j$ we get a morphism
$$
p_{ij}: G\times_SZ'_i@>{(id_G,s_{ij})}>> G\times_S X'_i @>m>> X'_i.
$$
Let $Z''_i$ denote the disjoint union of $j(i)$ copies of $Z'_i$. The sum
of the  morphisms $p_{ij}$  gives a finite and surjective $G$-morphism
$p_i\deq \sum_jp_{ij}$
and a commutative diagram
$$
\CD
G\times_SZ''_i@>{p_i}>> X'_i\\
@VVV @VVf'_iV\\
Z''_i@>>> Z'_i.
\endCD
$$
(3.14.1--3) show that it is sufficient to prove the following:

\proclaim{3.14.5 Lemma} (3.12) is true in case
$X=G\times_SZ$ and $f:X\to Z$ is the second projection.
\endproclaim

\demop
Since $G/S$ is affine, so is $G\times_SZ\to Z$. If $e:S\to G$
denotes the unit section then $Z\cong S\times_SZ\hookrightarrow X$, so $Z$ is
 of finite type  if $X$ is. These show (3.12.1) and (3.12.2).

$f:X\to Z$ is a group scheme over $Z$. Thus
$E\mapsto f^*E$ provides an equivalence of the category of (quasi)coherent
sheaves on $Z$ and the category  of (quasi)coherent
$G$-sheaves on $X$. Its inverse is given by
$f^*E\mapsto (f_*(f^*E))^G=E$.  Thus if $F$ is a coherent $G$-sheaf on $X$ then
$(f_*F)^G$ is a coherent sheaf on $Z$.\qed\enddemo

This completes the proof of (3.12).\qed\enddemo

\head 4. Quotients of Normal Spaces by Seshadri
\endhead

The aim of this section is to consider quotients of normal algebraic spaces
by proper actions of algebraic group schemes.  This situation has been
studied in \cite{Seshadri63,72} and his methods
 give an almost complete answer.

We need to apply his techniques in a more general setting. Therefore we
give precise statements.  The proofs require only slight technical
modifications. Some of these are worked out in \cite{Koll\'ar95a}.

\proclaim{4.1 Proposition} \cite{Seshadri72, Lemma 6.1} Let $G$ be an
algebraic group scheme,  smooth and of finite type over $S$.
Let $X$ be an algebraic  $G$-space and $n:\bar
X\to X$ the normalization. Assume that $n$ is finite.

 Then $\bar X$ has a unique structure of an algebraic  $G$-space such that $n$
is a
$G$-morphism.\qed\endproclaim

\demo{4.1.1 Remark} There are some natural group schemes over
$\z$ which are not smooth over $\z$. The simplest example is the group of
square roots of unity $\mu_2\deq \spec \z[x]/(x^2-1)$. This has two irreducible
components (corresponding to $+1$ and $-1$) which intersect in characteristic
2. Its normalization consists of two copies of $\spec \z$, and  the action of
$\mu_2$ on itself does not lift to an action on the normalization.
\enddemo

\proclaim{4.2  Theorem}  \cite{Seshadri72, Theorem 6.1} Let $X$ be an
excellent   algebraic space  over
$S$ and
$G$   an affine algebraic group scheme, smooth and of finite type over $S$. Let
$m_X:G\times X\to X$ be a proper
$G$-action on
$X$.

Then there is a Zariski locally trivial principal $G$-bundle $p:Y\to Z(Y)$
and a finite and dominant $G$-morphism $q:Y\to X$.

If $X$ is irreducible, then we can choose $p:Y\to Z(Y)$ and $q:Y\to X$ in such
a way that $Y$ is irreducible, the field extension $k(Y)\supset k(X)$ is normal
with Galois group
$\Gamma$ and the action of $\Gamma$  on $Y$ commutes with the $G$-action on
$Y$.\qed
\endproclaim

As a consequence we obtain the first quotient theorem:

\proclaim{4.3 Theorem} Let  $G$ be an  affine
algebraic group scheme, smooth and    of finite type over $S$ and $X$ a normal
 algebraic space  of finite type  over $S$. Let $m:G\times X\to X$
be a proper
$G$-action on $X$.

Then the geometric quotient $f:X\to X/G$ exists, and $X/G$ is a  normal
algebraic space    of finite type over $S$.
\endproclaim

\demop  We may assume that $X$ is irreducible. Choose $q:Y\to X$ and $p:Y\to
Z(Y)$ as in (4.2).  Set $Y'=Y/\Gamma$ and $Z'=Z(Y)/\Gamma$; these exist by
(4.4). $p$ descends to a morphism $p':Y'\to Z'$.  Since the $G$ and the
$\Gamma$ actions commute, the $G$-action on $Y$ descends to a $G$-action on
$Y'$.
{}From (3.11) we conclude that
$p'$ is a topological quotient of $Y'$ mod $G$.

Assume next that $k(X)$ has characteristic zero. Then $k(X)=k(Y)^{\Gamma}$,
hence  $Y'=Y$ and  $Z=Z'$ is a topological quotient of $Y$ mod $G$.
Thus a geometric quotient exists by (3.13).

If $k(X)$ has positive characteristic   then $Y'\to Y$ is a purely inseparable
$G$-morphism. A topological quotient of $Y'$ mod $G$ exists. By (6.7) this
implies that  a topological quotient of $Y$ mod $G$ also exists.\qed\enddemo

In the proof we used  that quotients of algebraic spaces by finite groups
exist. This result is due to Deligne (cf. \cite{Knutson71, p.185}), but I do
not know of any published proof. In the case that is used here,  this also
follows from (4.7).

\proclaim{4.4 Lemma}  Let $X$ be a  separated algebraic space
of finite type over $S$ and $G$ a finite
group acting on $X$. Then a geometric quotient $X/G$ exists and $X/G$ is also
of finite type over $S$.\qed
\endproclaim

It is not  more difficult  to understand quotients of normal algebraic spaces
by finite equivalence relations. Their existence is established through the
following two constructions.

\demo{4.5 Construction of Symmetric Products} Let $X$ be a separated
algebraic space. Our aim is to construct its symmetric products, denoted by
$S^nX$. $X$ can be given by an \'etale equivalence relation
$R\rightrightarrows U$ where $U$ is affine. Since
$X$ is separated, the resulting morphism $\phi:R\to U\times U$ is a closed
immersion. We can write $R=\Delta\cup R'$ (disjoint union) where
$\phi(\Delta)$ is the diagonal.

The induced morphism $f^n:U^n\to X^n$ is \'etale. The symmetric group $S_n$
acts on both spaces and  $f^n$ is $S_n$-equivariant.  Unfortunately, it is
not fixed-point reflecting. More precisely, it is not fixed-point reflecting
at $(x_1,\dots,x_n)$ iff there is a pair $i<j$ such that $f(x_i)=f(x_j)$ and
$x_i\neq x_j$.
Set
$$
W_n\deq U^n\setminus \bigcup_{i<j} \pi_{ij}^{-1}(\phi(R')),
$$
where $\pi_{ij}:U^n\to U^2$ is the projection to the product of the
$i^{\text{th}}$ and
$j^{\text{th}}$ factors.    $f^n:W_n\to  X^n$ is \'etale, surjective and
fixed-point reflecting.  $W_n$ is quasi affine, hence the geometric quotient
$W_n/S_n$  exists and is  also quasi affine (cf. \cite{Mumford68, Sec. 12}).
Thus
$S^nX$ is given by the equivalence relation
$$
(W_n\times_{X^n}W_n)/S_n \rightrightarrows W_n/S_n.
$$
\enddemo

\proclaim{4.6 Lemma}  Let $X,Y$ be separated algebraic spaces, $X$ normal
and irreducible. Let $Z\subset X\times Y$ be a closed subspace such that the
projection
$p_X:Z\to X$ is finite and surjective on every irreducible component of $Z$.
Let $d=\deg (Z/X)$. Then there is a natural morphism
$\phi_Z:X\to S^d Y$.
\endproclaim

\demop Let $X^0\subset X$ be a dense open set such that $p_X|Z$ is flat over
$X^0$. Then
$$
x\mapsto p_Y[p_X^{-1}(x)]\in  S^dY
$$
is a well defined morphism. By looking at the closure of the graph it is easy
to see that it extends to $\phi_Z:X\to S^d Y$ such that
$\supp \phi_Z(x)=\supp p_Y[p_X^{-1}(x)]$ for every $x\in X$. (In general no
such extension exists if
$X$ is not normal.)\qed\enddemo

\proclaim{4.7 Corollary}  Let $X$ be a normal algebraic space. Let $Z\subset
X\times X$ be a closed subspace such that the first projection
$p_X:Z\to X$ is finite and surjective on every irreducible component of $Z$.
  Assume that $Z$ defines an equivalence relation on
$X$. (For the present purposes all we need is that it is an equivalence
relation on closed points in algebraically closed fields.) Then there is a
 morphism
$\phi_Z:X\to X/Z$ whose fibers are precisely the $Z$-equivalence classes.
\endproclaim

\demop Set  $d=\deg (Z/X)$ and let $X/Z$ be the normalization of the image of
$\phi_Z:X\to S^d X$. \qed\enddemo

\demo{4.7.1 Remark} The above construction gives the correct fiber set
theoretically, but not necessarily scheme theoretically, especially in
positive characteristic. Thus (4.7) implies that if $X$ is a normal algebraic
space and $G$ a finite and flat group scheme acting on $X$ then a topological
quotient exists. Thus  by (3.13) a geometric quotient also exists.
\enddemo

\head 5. Weakly Normal Spaces
\endhead

The aim of this section is to consider the quotient problem for weakly normal
algebraic spaces. The notion of weak normality and of universal
homeomorphisms are central to the construction of quotients.  We start by
recalling their definitions and  basic properties.

\demo{5.1 Definition} We say that a morphism of schemes $g:U\to V$ is a {\it
universal homeomorphism} if  it is a homeomorphism and for every $W\to V$ the
induced morphism $U\times_VW\to W$ is again a homeomorphism.

This notion extends to morphisms of algebraic spaces the usual way
\cite{Knutson71, II.3}.
\enddemo

This concept is especially useful if $g$ is finite. In this case we have the
following   characterization (cf.
\cite{EGA71, I.3.7--8}).

\proclaim{5.2 Lemma}  Let $g:U\to V$ be a finite and surjective morphism of
algebraic spaces. The following are equivalent

(5.2.1) $g$ is  a universal homeomorphism.

(5.2.2) $g$ is surjective and universally injective.

(5.2.3) For every $v\in V$ the fiber  $g^{-1}(v)$ has a single
point $v'$ and
$k(v')$ is a purely inseparable field extension of   $k(v)$.
(Thus $k(v')=k(v)$ if $\chr k(v)=0$.)\qed
\endproclaim

\demo{5.2.4 Remark} The terminology in English does not seem to be uniform.
The french ``radiciel" \cite{EGA71, I.3.7} did not catch on. Some authors use
the name ``purely inseparable morphism" to refer to a morphism satisfying the
equivalent properties in (5.2).
\enddemo

One of the most important properties of these morphisms is the following:

\proclaim{5.3 Theorem} \cite{SGA1, IX.4.10}  Let $g:U\to V$ be a finite and
universal homeomorphism. Then $\ast\mapsto \ast\times_VU$ induces an
equivalence between the categories
$$
(\text{\'etale morphisms: $\ast\to V$})
@>{\ast\mapsto \ast_VU}>>
(\text{\'etale morphisms: $\ast\to U$}).\qed
$$
\endproclaim

The notion of weak normality
was developed in a series of papers
\cite{Andreotti-Noguet67; Andreotti-Bombieri69;
Traverso70}. We recall the definitions and  a few of the basic properties.
A short outline of the theory can also be found in \cite{Koll\'ar95b, I.7}.

\demo{5.4 Definition}
The {\it weak normalization}  of an algebraic space $X$ is an algebraic space
$ X'$  together with a finite and surjective morphism $m: X'\to X$ such that

(5.4.1) for every $x\in X$ there is a unique $ x'\in X'$ such that
$m( x')=x$;

(5.4.2) for every  $x\in X$ the field extension $k(x')\supset
k(x)$ is purely inseparable and an isomorphism if $x\in X$ is a generic point;

(5.4.3) if $g:Y\to X$ satisfies the above two properties then there is a
unique  factorization
$$
m: X'@>>> Y@>g>> X.
$$

 An algebraic space $X$ is called  {\it weakly normal}   if it is its own
 weak  normalization.
\enddemo

\proclaim{5.5 Proposition} \cite{Andreotti-Noguet67, Andreotti-Bombieri69,
Traverso70}
 Let $X$ be an algebraic space  such that its
 normalization $\bar X\to X$ is finite over $X$.
Then:

(5.5.1) $X$ has a   weak normalization
 $n:X^{\text{wn}}\to X$. $n$ is a finite and  universal homeomorphism.

(5.5.2) There is a factorization  $\bar X\to X^{\text{wn}}\to X$.

(5.5.3) The conductor ideal (5.7)  $I_{\bar X/X^{\text{wn}}}\subset \o_{\bar
X}$ is
its own radical.

(5.5.4) Let $f:Y\to X$ be a morphism such that $X$ is weakly normal at the
image of every generic point of $Y$. Then $f$ lifts to a morphism
$f^{\text{wn}}:Y^{\text{wn}}\to X^{\text{wn}}$.

(5.5.5) If $f:Y\to X$ is smooth and $X$ is weakly normal then so is $Y$.
\endproclaim

\demop  For the reader's convenience, let us recall the construction of $
X^{\text{wn}}$ for schemes.

For a point $x\in X$ let $\bar x=n^{-1}(x) \subset \bar X$
with reduced scheme structure.   $k(\bar x)$ is the direct sum of the residue
fields of the points of $\bar x$, $i_x:k(x)\to k(\bar x)$   the
natural injection. For an open affine $U\subset X$ let
$$
\o^{\text{wn}}(U)\deq \{f\in \o_{n^{-1}(U)}| f(\bar x)^{q}\in
\im i_x\ \forall x\in U,
\text{ where $q$ is a power of $\chr k(x)$.}\}
$$
This defines a sheaf of rings  $\o^{\text{wn}}$ over $X$. Set
$ X^{\text{wn}}=\spec_X\o^{\text{wn}}$.

The properties (5.5.1--4) for schemes are clear from this construction.
(5.5.5) is in \cite{Koll\'ar95b, I.7.2.6}.

Let $u:U\to V$ be a surjective \'etale morphism of schemes. By (5.5.5) $U$
is weakly normal iff $V$ is.  Thus all the
above assertions descend from the category of schemes to the category of
algebraic spaces (cf. \cite{Knutson71, p.106}).\qed\enddemo

\proclaim{5.6 Corollary} Let $G$ be an
algebraic group scheme,  smooth and of finite type over $S$.
Let $X$ be an algebraic  $G$-space and $n:
X^{\text{wn}}\to X$ the weak normalization. Assume that $n$ is finite.

 Then $X^{\text{wn}}$ has a unique structure of an algebraic  $G$-space
such that $n$ is a
$G$-morphism.\endproclaim

\demop The $G$-action on $X$ is given by a morphism $m:G\times X\to X$.
Since $G$ is smooth over $S$, (5.5.5) implies that
$G\times (X^{\text{wn}})=(G\times X)^{\text{wn}}$. Thus by (5.5.4) $m$ lifts
to a morphism
$$
m^{\text{wn}}: G\times X^{\text{wn}}\to X^{\text{wn}}.
$$
It is clear that $m^{\text{wn}}$ defines a group action.\qed\enddemo

Our aim is to construct geometric quotients of weakly normal $G$-spaces by
comparing the situation to the normalization. Let $\bar X\to X$ be the
normalization. Assume that  the $G$-action lifts to an action on $\bar X$ and
that a geometric quotient $\bar X/G$ exists. Thus we have a diagram
$$
\CD
\bar X@>>> X\\
@VVV @.\\
\bar X/G @.
\endCD
$$
and $X/G$ would complete it to a commutative square. The construction of
$X/G$ becomes a special case of constructing push-outs (8.1).

Before   stating the main result of this section, we need to introduce
further notation:

\demo{5.7 Notation}  Let $G/S$ be an algebraic group scheme, flat over $S$.
Let
$X/S$ be a weakly normal algebraic space  over $S$ with normalization
$p:\bar X\to X$ which is assumed to be finite over $X$.
Assume that we are given
  $G$-actions on $\bar X$ and on $X$ such that $n$ is a $G$-morphism.

Let $I_X\subset \o_X$ resp. $I_{\bar X}\subset \o_{\bar X}$ denote the
{\it conductor ideals}. That is,
$$
I_{\bar X}=I_X=\hom_X(p_*\o_{\bar X}, \o_X).
$$
This definition shows that $I_X$ and $I_{\bar X}$ are $G$-sheaves.
The corresponding    subschemes are denoted by $C_X\subset X$
resp.
$C_{\bar X}\subset \bar X$.
By (5.5.3)   $I_X$ and $I_{\bar X}$ are their own radicals,
so the conductor subschemes are reduced. The $G$ action restricts to
$G$-actions on
$C_X$ resp.
$C_{\bar X}$.
$p$ restricts to a finite and surjective $G$-morphism
$p:C_{\bar X}\to C_X$.

Assume that a geometric quotient $\bar f:\bar X\to Z(\bar X)$ exists.
$C_{\bar X}\subset \bar X$ is $G$-invariant and closed, thus there is a
closed and reduced subscheme  $j:Z(\bar C)\hookrightarrow  Z(\bar X)$  such
that
$\bar f^{-1}(Z_{\bar C})=C_{\bar X}$ (as sets, but not necessarily as
schemes).

By (3.11.3)
$C_{\bar X}\to Z(\bar C)$ is a topological quotient of $C_{\bar X}$
by
$G$. Therefore, by (3.13) a geometric quotient  $C_{\bar X}\to Z(C_{\bar
X})$ exists and the induced morphism $q_1: Z(C_{\bar X}) \to   Z(\bar C)$
is a finite and universal homeomorphism.

Assume furthermore that
a geometric quotient  $C_X\to Z(C_X)$  also exists. Since $p:C_{\bar X}\to
C_X$ is a $G$-morphism, we obtain an induced morphism $q_2: Z(C_{\bar X}) \to
Z(C_X)$.
\enddemo

\proclaim{5.8 Theorem}   Let $G$ be an algebraic group scheme, flat and of
finite type over $S$. Let $X$ be a  weakly normal
algebraic $G$-space of finite type over $S$ such that the $G$-action lifts to a
$G$-action on the normalization $p:\bar X\to X$.  Assume that:

(5.8.1) The  geometric quotients
 $\bar f:\bar X\to Z(\bar X)$ and $C_X\to Z(C_X)$ exist.

(5.8.2) The push-out diagram
$$
\CD
Z(C_{\bar X})@>q_1>>  Z(\bar C)\\
@Vq_2VV @.\\
Z(C_X) @.
\endCD\tag 5.8.2.1
$$
can be completed to
$$
\CD
Z(C_{\bar X})@>q_1>>  Z(\bar C)\\
@Vq_2VV @VVs_2V\\
Z(C_X) @>s_1>> W,
\endCD\tag 5.8.2.2
$$
such that $s_2$ is   finite and  $s_1$ is a
finite and universal homeomorphism.
\smallpagebreak

Then a geometric quotient of $X$ mod $G$ exists.
\endproclaim

\demop Assume that $ s_2$ is finite. By (8.6)
there is an algebraic space $Z$ and  finite morphisms
$q:Z(\bar X)\twoheadrightarrow Z$  and $j_W: W\hookrightarrow Z$ which fit into
the following commutative diagram
$$
\CD
Z(\bar C)@>j>> Z(\bar X)\\
@Vs_2VV @VVqV\\
W@>j_W>> Z.
\endCD
$$
The next two claims imply (5.8):

\proclaim{5.8.3 Claim}  (5.8.3.1) $\bar f: \bar X\to Z(\bar X)$ descends to
a
$G$-morphism
$f:X\to Z$.

(5.8.3.2) $f:X\to Z$ is a topological quotient of $X$ by $G$.
\endproclaim

\demop  The first assertion means that the diagram
$$
\CD
\bar X@>p>> X\\
@V\bar fVV @.\\
Z(\bar X)@>q>> Z
\endCD
$$
can be completed to a push-out diagram by  $f:X\to Z$. At the set theoretic
level this defines $f$ uniquely as $f=q\circ \bar f\circ p^{-1}$.

 We construct $f$ by constructing its graph. The morphism
$p\times (q\circ \bar f):\bar X@>>> X\times Z$ is finite; let $Y\subset
 X\times Z$ denote its image (with reduced scheme structure).

$Y$ is the graph of the (a priori multi-valued) map
$q\circ \bar f\circ p^{-1}$.
If  the projection $h:Y\to X$ is an isomorphism, then it defines a morphism
$f:X\to Z$.
$X$ is weakly normal and $Y$ is dominated by the normalization $\bar
X$, thus by (5.5) it is sufficient to prove that $h:Y\to X$ is universally
injective.

Given $x\in X$ let $\red p^{-1}(x)=\bar x_1\cup\dots \cup\bar x_s$ (where $s$
depends on $x$).   Then
$$
\red h^{-1}(x)=\bigcup_i(x,(q\circ \bar f)(\bar x_i)).
$$
 If $x\not\in C_X$ then $\bar X\to X$ is an
isomorphism over $x$, thus $Y\to X$ is also an isomorphism
over $x$.

Next assume that $x\in C_X$, hence $\red p^{-1}(x)\subset C_{\bar X}$.
$q\circ \bar f:C_{\bar X}\to Z$ sits in a commutative diagram
$$
\CD
C_{\bar X}@>\bar f>>Z(C_{\bar X})@>q_1>>  Z(\bar C)@>j>>Z(\bar X)\\
@VpVV @Vq_2VV @VVs_2V @VqVV\\
C_X@>>> Z(C_X) @>s_1>> W @>j_W>> Z.
\endCD\tag 5.8.3.3
$$
This shows that $(q\circ \bar f)(\bar x_i)$ is independent of $i$ and the
residue field of $Z$ at $(q\circ \bar f)(\bar x_i)$ is a subfield of the
residue field of $X$ at $x$. Therefore
$h:Y\to X$ is  universally injective, hence an isomorphism. This shows
(5.8.3.1).

The morphisms $p:\bar X\to X$ and $q: Z(\bar X)\to Z$ are finite, thus $Z$ is
an approximate quotient of $X$ mod $G$ by (3.11).

In order to see that it is  a topological quotient we need to check that for
every algebraically closed field $K$ (2.7.2.3) holds. By construction
$$
\align
X(K)/G(K)&=(X-C_X)(K)/G(K)\cup C_X(K)/G(K)\\
&=(\bar X-C_{\bar X})(K)/G(K)\cup
C_X(K)/G(K) \\
&=(Z(\bar X)-Z(\bar C))(K)\cup Z(C_X)(K)\\
&=(Z-W)(K)\cup W(K)=Z(K).
\endalign
$$
All the unions are disjoint, thus $X(K)/G(K)\to Z(K)$ is an isomorphism of
sets.
\qed\enddemo
\enddemo

\head 6. Quotients in Positive Characteristic
\endhead

In this section, $k$  denotes a field of positive characteristic $p$.
The main result of the section is the following:

\proclaim{6.1 Theorem} Let  $G$ be an affine
algebraic group scheme   of finite type over $k$ and $X$ an algebraic
space  of finite type over $k$. Let $m:G\times X\to X$ be a proper $G$-action
on $X$.

Then a geometric quotient $f:X\to X/G$ exists, and $X/G$ is an
algebraic
space  of finite type over $k$.
\endproclaim

We already proved this result in case $X$ is   normal (4.2). In this
section we reduce the general case to the   normal case.

The crucial feature of positive characteristic algebraic spaces that we
use is the Frobenius morphism:

\demo{6.2 Definition of the Geometric Frobenius Morphism \cite{SGA5, XIV}}

Let $S$ be an $\f_p$-scheme.
Fix $q=p^r$ for some natural number $r$. Then $a\mapsto a^q$ is an
endomorphism of $\o_S$. This can be extended to polynomials
by the formula
$$
f=\sum a_Ix^I\mapsto f^{(q)}\deq \sum a^q_Ix^I.
$$

Let $U=\spec R$ be an affine scheme over $S$. Write
$R=\o_S[x_1,\dots,x_m]/(f_1,\dots,f_n)$ and set
$$
R^{(q)}\deq \o_S[x^{(q)}_1,\dots,x^{(q)}_m]/(f^{(q)}_1,\dots,f^{(q)}_n)
\qtq{and}  U^{(q)}\deq
\spec R^{(q)},
$$
where the $x^{(q)}_i$ are new variables.
There are natural morphisms
$$
F^q:U\to U^{(q)}\qtq{and} (F^q)^*:R^{(q)}\to R\qtq{given by}
(F^q)^*(x^{(q)}_i)=x_i^q.
$$

It is easy to see that these are independent of the choices made. Thus $F^q$
gives a functor from algebraic spaces over $S$  to algebraic spaces over
 $S$.

If $F^q:\spec S\to \spec S$ denotes the $q^{\text{th}}$-power morphism, then
one can define $X^{(q)}$ intrinsically as
$$
X^{(q)}=X\times_{\spec S,F^q}\spec S.
$$
\enddemo

Since $X^{(q)}$ is obtained by base change,  it commutes with many other
constructions. In particular, we obtain the following:

\proclaim{6.3 Lemma} Let $m:G\times X\to X$ be a $G$-action.

(6.3.1) $m^{(q)}:G^{(q)}\times X^{(q)}\to X^{(q)}$ defines a $G^{(q)}$-action
on $X^{(q)}$.

(6.3.2) $m$ is proper iff $m^{(q)}$ is.

(6.3.3) $f:X\to Z$ is a topological quotient of $X$ by $G$ iff
$f^{(q)}:X^{(q)}\to Z^{(q)}$ is a  topological quotient of $X^{(q)}$ by
$G^{(q)}$.\qed
\endproclaim

Frobenius morphisms give the most interesting examples of
universal homeomorphisms:

\proclaim{6.4 Lemma}  Let $X$ be an algebraic space of finite type over
$S$. Then
$F^q:X\to X^{(q)}$ is a finite and universal homeomorphism.\qed
\endproclaim

\demo{6.5 Remark} There are many algebraic spaces which are not of finite
type and such  that $F^q:X\to X^{(q)}$ is  finite. By a result of Kunz (see
\cite{Matsumura80, p.302}) such algebraic spaces are excellent.
The results of this section  probably  remain valid for algebraic spaces
satisfying this property. The following example of \cite{Nagata69} shows that
this assumption is necessary.
\enddemo

\demo{6.5.1 Example} Let $k$ be a field of characteristic
$p>0$  and $K\deq k(x_1,x_2,\dots)$, where the
$x_i$ are algebraically independent over $k$.  Let
$$
D\deq \sum_ix_i^{1+p}\frac{\partial}{\partial x_i}\qtq{be a derivation of $K$.}
$$
Let $F\deq\{f\in K\vert D(f)=0\}$ be the subfield of constants.
Set
$$
R=K+\epsilon K\qtq{where $\epsilon^2=0$, and} \sigma: f+\epsilon g\mapsto
f+\epsilon (g+D(f)).
$$
$R$ is a local Artin ring.
It is easy to check that $\sigma$ is an automorphism of $R$ of order $p$.
The fixed ring is $R^{\sigma}=F+\epsilon K$. Its maximal ideal is
$m\deq (\epsilon K)$ and  generating sets of $m$ correspond to  bases
of $K$ as an $F$-vectorspace. (6.5.2) shows that $m$ is not
finitely generated.

\proclaim{6.5.2 Claim} Notation as above. The $x_i$ are linearly independent
over
$F$.
\endproclaim

\demop Assume that we have a relation $\sum_{i\in I} f_ix_i=0$
where $|I|$ is the smallest possible. Apply $D$ to
get that
$$
0=\sum_{i\in I} f_iD(x_i)=\sum_{i\in I} (x_i^pf_i)x_i.
$$
$x_i^p\in F$, thus
we get another linear dependence relation, a contradiction.
\qed\enddemo\enddemo

For us the most important feature of the Frobenius morphism is the
following universal property:

\proclaim{6.6 Proposition}  Let $X,Y$ be  algebraic spaces over  $S$  and
$g:X\to Y$ a finite  and universal homeomorphism. Then there is a $q=p^r$
such that   $F^q$ can be factored as
$$
F^q:X@>g>>Y @>\bar g>> X^{(q)}.
$$

For large $q$ the factorization is canonical in the sense that if
$$
\CD
X_1@>g_1>> Y_1\\
@VfVV @VhVV\\
X_2@>g_2>> Y_2\\
\endCD
$$
is a commutative diagram where the  $g_i$ are
finite and universal homeomorphisms, then the factorization gives a
commutative diagram
$$
\CD
X_1@>g_1>> Y_1@>\bar g_1>> X_1^{(q)}\\
@VfVV @VhVV @Vf^{(q)}VV \\
X_2@>g_2>> Y_2@>\bar g_2>> X_2^{(q)}.
\endCD
$$
\endproclaim

\demop Since  the factorization is canonical, it is sufficient
to construct it in case $X$ and $Y$ are affine schemes over an affine scheme
$\spec C$.

The latter can be formulated in terms of a ring homomorphism $g^*:A\to B$ where
$A$ and $B$ are $C$-algebras.
We can decompose $g^*$ into $A\twoheadrightarrow B_1$ and $B_1\hookrightarrow
B$. We deal with them separately.

First consider $B_1\subset B$. We need to show that there is a $q$ such that
$B^q\subset B_1$, where $B^q$ denotes the $C$-algebra generated by the $q$-th
powers of all elements. The  proof is by Noetherian induction.

If $B_1$ is Artin, then the assertion is clear. Applying this over the
generic points we obtain that  $B_1\subset B_1B^q$ is an isomorphism at all
generic points for
$q\gg 1$. Let $I<B_1$ denote the conductor of this extension.

By induction we know that there is a $q'$ such that
$(B_1B^q/I)^{q'}\subset B_1/I$. Thus we get
$$
B^{(qq')}\to B^{qq'}\subset (B_1B^q)^{q'}\subset B_1.
$$

Next consider $A\twoheadrightarrow  B_1$.  The kernel
 is a nilpotent ideal $I<A$, say $I^m=0$. Choose $q'$ such that $q'\geq
m$. For $b_1\in B_1$ let $b_1'\in A$ be any preimage. Then
$(b'_1)^{q'}$ depends only on $b_1$. The map
$$
b_1\mapsto (b'_1)^{q'}\qtq{defines a factorization}
B_1^{(q')}\to A\to B_1.
$$
Combining the map  $B^{(q)}\to B_1$ with  $B_1^{(q')}\to A$ we
obtain  $B^{(qq')}\to A$.\qed
\enddemo

The main technical result of the section is the following:

\proclaim{6.7 Theorem}  Let $G,H$ be algebraic group schemes of finite type
over $S$ and $\rho:H\to G$ a group homomorphism which is a finite and universal
homeomorphism. Let $Y,X$ be algebraic spaces of finite type
over $S$ and  $f:Y\to X$ a finite and universal homeomorphism.

Let $m_Y:H\times Y@>>> Y$ and $m_X:G\times X@>>> X$ be proper actions and
assume that $f$ is a $\rho$-morphism, that is, the following diagram is
commutative
$$
\CD
H\times Y@>{\rho\times f}>> G\times X\\
@Vm_YVV @Vm_XVV\\
Y@>f>> X.
\endCD\tag 6.7.1
$$
Then the following are equivalent:

(6.7.2) The geometric quotient of $X$ mod $G$ exists.

(6.7.3) The geometric quotient of $Y$ mod $H$ exists.
\endproclaim

\demop  Let
$p_X:X\to X/G$ be the geometric quotient. Then
$p_X\circ f:Y\to X/G$ is a topological quotient by (3.11), hence a
geometric quotient of $Y$ mod $H$ exists by (3.13).

Conversely,  let $p_Y:Y\to Y/H$ be the geometric quotient. By
(6.3.3)
$p_Y^{(q)}:Y^{(q)}\to (Y/H)^{(q)}$ is a topological quotient of $Y^{(q)}$ mod
$H^{(q)}$. By (6.6) the diagram (6.7.1) can be completed to
$$
\CD
H\times Y@>{\rho\times f}>> G\times X
@>{\bar{\rho}\times \bar f}>> H^{(q)}\times Y^{(q)}\\
@Vm_YVV @Vm_XVV @Vm_Y^{(q)}VV \\
Y@>f>> X@>\bar f>> Y^{(q)}.
\endCD\tag 6.7.4
$$
$p_Y^{(q)}\circ \bar f:X\to (Y/H)^{(q)}$ is a topological quotient as before,
hence a  geometric quotient of $X$ mod $G$ exists by (3.13). \qed
\enddemo

\proclaim{6.8 Corollary}  Let $G$ be an affine algebraic group scheme of
finite type over $k$ and $X$ an algebraic $G$-space of finite type
over $k$. Assume that the $G$-action on $X$ is proper.

Then, a geometric quotient of $X$ mod $G$ exists iff
a geometric quotient of $X$ mod $\red G$ exists.
\endproclaim

\demop Apply (6.7) with $Y=X$ and $H=\red G$.\qed\enddemo

\proclaim{6.9 Corollary}  Let $G$ be a smooth, affine  algebraic group
scheme of finite type over $k$ and $X$ an algebraic $G$-space of finite type
over $k$. Assume that the $G$-action on $X$ is proper.

Then, a geometric quotient of $X$ mod $G$ exists iff
a geometric quotient of $(\red X)^{\text{wn}}$ mod $G$ exists.
\endproclaim

\demop Apply (6.7) with $Y=(\red X)^{\text{wn}}$ and $H=G$.
The $G$ action on $X$ gives a $G$-action on $\red X$ and this in turn gives a
$G$-action on $(\red X)^{\text{wn}}$ by (5.6). (Both times we used that $G$ is
smooth.)
\qed\enddemo

\proclaim{6.10 Corollary}
 Let $G$ be a smooth, affine algebraic group scheme of finite
type over $k$ and $X$ a weakly normal algebraic $G$-space of finite type
over $k$. Assume that the $G$-action on $X$ is proper.
Let $p:\bar X\to X$ denote the normalization.

Then, the geometric quotient of $X/G$ exists if
 $\bar X/G$ exists.
\endproclaim

\demop Assume that  $\bar f:\bar X\to \bar X/G$ exists.  We use the notation
of (5.7).

$Z(C_{\bar X})\to Z(\bar C)$ is a finite and universal homeomorphism by
(3.11.3).  Hence by (8.4)  a push-out diagram as in (5.8.2.2) exists. Thus a
geometric quotient of $X$ mod $G$ exists by (5.8).\qed\enddemo

\demo{6.11 Proof of (6.1)}  Start with arbitrary $G$ and $X$. We would like
to construct a geometric quotient of $X$ by $G$. By (6.3.3) it is sufficient
to show that the geometric quotient   $X^{(q)}/G^{(q)}$ exists for some $q$.
If
$q\gg 1$ then $\red G^{(q)}$ is smooth over $k$ by \cite{Weil62, p.25}. Thus by
(6.8) it is sufficient to show the existence of a geometric quotient
in case $G$ is smooth over $k$.

By (6.9) this can be further reduced to the case when in addition $X$ is
weakly normal, and by (6.10) to the case when $X$ is normal. The latter types
of actions have already been dealt with in (4.2).\qed
\enddemo

The methods of this section also apply to group schemes in positive
characteristic. The only part that fails is the reduction
to smooth group schemes in (6.11). We obtain the following result:

\proclaim{6.12 Theorem} Let $S$ be an excellent scheme of characteristic $p$.
Let  $G$ be a smooth, affine
algebraic group scheme   of finite type over $S$ and $X$ an algebraic
space  of finite type over $S$. Let $m:G\times X\to X$ be a proper $G$-action
on $X$.

Then the geometric quotient $f:X\to X/G$ exists, and $X/G$ is an
algebraic
space  of finite type over $S$.\qed
\endproclaim

\head 7. Quotients by Reductive Group Schemes
\endhead

\demo{7.1 Definition} A  group scheme $G/S$ is called reductive if $G$ is
smooth over $S$ and for every field $\spec K\to S$ the fiber product $G_K\deq
G\times_S\spec K$ is a connected and reductive algebraic group over $K$.

See \cite{Demazure67} for a detailed study of the structure of reductive group
schemes.
\enddemo

\demo{7.2 Examples} The simplest and for moduli problems the most important
  reductive group schemes are  $SL(n,\z)$ and  $PGL(n,\z)$, as  group schemes
over $\z$. For a field $K$ the corresponding algebraic groups are $SL(n,K)$
resp. $PGL(n,K)$.
\enddemo

\cite{Seshadri77}  studied the invariant theory of reductive group schemes.
Two of his results are crucial for our purposes. (As always, we assume that
$S$ is excellent.)

\proclaim{7.3 Theorem} \cite{Seshadri77} Let  $G/S$ be a
reductive group scheme.  Let $X/S$ be an affine scheme of finite type over
$S$ and $m:G\times X\to X$ a $G$-action. Assume for simplicity that the
action is proper. Then the following hold.

(7.3.1) The geometric quotient $f:X\to X/G$  exists and $X/G$ is
an affine scheme of finite type over
$S$.

(7.3.2) If $S$ is a scheme over $\spec \q$ then  $f:X\to X/G$ is a
universal geometric quotient.  That is, if
$g:Y\to X$ is a $G$-morphism which is a
 closed immersion then  the induced morphism  $f/G:Y/G\to X/G$ is  also a
closed immersion.
\endproclaim

\demop  The first assertion is \cite{Seshadri77, p.263}. The second part is
not mentioned explicitly, so I outline how it follows from the results of the
article.

Pushing forward to the base scheme $S$, (7.3.2)  is equivalent to
the following assertion:

\proclaim{7.3.3 Claim}  Let $S$ be a scheme over $\spec \q$ and  $G/S$  a
reductive group scheme. Let $0\to
F_1\to F_2\to F_3\to 0$ be an exact sequence of quasi coherent $G$-sheaves on
$S$.  Then the sequence
$$
0\to  F_1^G\to  F_2^G@>v>>  F_3^G\to 0
$$
is exact.
\endproclaim

\demop By (3.14.3.3) the above sequence is exact, except possibly at $F_3^G$.
Thus we need to establish that $v$ is surjective. By
\cite{Seshadri77, p.236} any quasi coherent $G$-sheaf is the union of  its
coherent $G$-subsheaves, thus it is sufficient to prove (7.3.3) for $F_3$
coherent. Then $F_3^G$ is also coherent, hence, by the Nakayama lemma, $v$ is
surjective if it is surjective over closed points. This
reduces (7.3.3) to the case when $F_3$ is a vector space over $k(s)$ for some
closed point $s\in S$.

A similar problem, with $S$ a scheme over $\spec \z$, is considered in\
\cite{Seshadri77, p.244}, where it is proved that it reduces to the problem
over $S=\spec\z$. The same argument shows that if $S$ is over $\spec \q$, then
(7.3.3) reduces to the special case where $S=\spec \q$.

Over  a field of characteristic zero any representation of a reductive group
is completely reducible, thus (7.3.3) holds for $S=\spec \q$. \qed\enddemo
\enddemo

\proclaim{7.4 Corollary} Let  $G/S$ be a
reductive group scheme and $S$  a scheme over $\spec \q$. Let $g:Y\to X$ be a
$G$-morphism of algebraic spaces over $S$ which is a closed immersion.
Assume that the  $G$-action on $X$ is proper and the geometric quotient $f:X\to
X/G$  exists.

Then the  geometric quotient $f:Y\to Y/G$  also exists
and the induced morphism  $f/G:Y/G\to X/G$ is  a closed immersion.
\endproclaim

\demop The question is local in the \'etale topology of $X/G$, thus we may
assume that $X/G$ is affine. By (3.12.1) then $X$ and $Y$ are both affine.
This is (7.3.2).\qed\enddemo

The main result of this section is the following:

\proclaim{7.5 Theorem} Let  $G$ be a reductive
algebraic group scheme    over $S$ and $X$ an algebraic
space  of finite type over $S$. Let $m:G\times X\to X$ be a proper $G$-action
on $X$.

Then the geometric quotient $f:X\to X/G$  exists  and $X/G$ is an
algebraic
space  of finite type over $S$.
\endproclaim

\demop We already proved this result in case $X$ is   normal (4.2). Here we
reduce the general case to the   normal case.

By Noetherian induction we may assume that the geometric quotient exists for
any $G$-invariant subscheme $Y\subset X$ (assuming $Y\neq X$).

First we consider
the weakly normal case, and then the general case.
\enddemo

\demo{7.6 Proof of (7.5) in case $X$ is weakly normal} We use the notation
 of (5.8).  We neeed to check that the two conditions (5.8.1--2) are
satisfied.

The geometric quotient  $\bar X/G$ exists by (4.2) and
  $C_X/G$ exists by the inductive assumption.
Thus (5.8.1) is satisfied.

Next look at the diagram  (5.8.2.1)
$$
\CD
Z(C_{\bar X})@>q_1>>  Z(\bar C)\\
@Vq_2VV @.\\
Z(C_X) @.
\endCD
$$
$C_{\bar X}\to \bar X$  is a closed immersion, thus the induced morphism
$Z(C_{\bar X})\to Z(\bar X)$ is a closed immersion at all points of
characteristic zero by (7.4). $Z(\bar C)$ is excatly the image of this
morphism, thus
$Z(C_{\bar X})@>q_1>>  Z(\bar C)$ is a  finite and  universal homeomorphism
which is an isomorphism  at all points where the characteristic of the residue
field is zero.

By (8.9) the above diagram can be completed to
$$
\CD
Z(C_{\bar X})@>q_1>>  Z(\bar C)\\
@Vq_2VV @VVs_2V\\
Z(C_X) @>s_1>> W.
\endCD
$$
Therefore, by (5.8),     $X/G$ exists.\qed\enddemo

\demo{7.7 Proof of (7.5) for arbitrary $X$}
Let $n:X^{\text{wn}}\to X$ be the weak normalization. By (5.5), $n$ is a
finite and  universal homeomorphism.

Let $f':X^{\text{wn}}\to X^{\text{wn}}/G$ be the geometric quotient. Let $U'$
be an affine scheme and $u:U'\to X^{\text{wn}}/G$ a surjective \'etale
morphism. Set
$$
V'=X^{\text{wn}}\times_{X^{\text{wn}}/G}U'.
$$
$f'$ is affine by (3.12.1), thus $V'$ is also affine. The first projection
$v':V'\to X^{\text{wn}}$ is an \'etale and fixed-point
reflecting $G$-morphism, since it is the pull back of the
\'etale morphism
$u$ (2.13).

By (7.8) the  fixed-point
reflecting, surjective and \'etale  $G$-morphism  $V'\to X^{\text{wn}}$
descends to a fixed-point reflecting, surjective and \'etale  $G$-morphism
$V\to X$.  $V'\to V$ is finite, hence $V$ is also affine.
The geometric quotient $V/G$ exists by (7.3.1). Thus by (2.17)
the geometric quotient $X/G$ also exists.\qed\enddemo

\proclaim{7.8 Lemma}  Let $G/S$ be a group scheme and
 $g:X'\to X$  a $G$-morphism which is a finite and
universal homeomorphism. Then $V\mapsto V\times_XX'$ induces an
equivalence between the categories
$$
(\text{\'etale $G$-morphisms: $\ast\to X$})
@>{V\mapsto V\times_XX'}>>
(\text{\'etale $G$-morphisms: $\ast\to X'$}).
$$
Moreover, $V\to X$ is fixed-point
reflecting iff $V\times_XX'\to X'$ is.
\endproclaim

\demop Let $V'\to X'$ be an \'etale $G$-morphism. By (5.3) there is an
\'etale morphism $V\to X$ such that $V'=V\times_XX'$. First we need to prove
that
$G$ acts on $V$ in the expected way. A $G$-action on $V'$ is given by a
morphism
$m':G\times V'\to V'$. We can encode this information as an isomorphism
$$
M': G\times V'@>{(id_G,m')}>> G\times V'.
$$
$G\times V'\to G\times X'$ is \'etale, thus by (5.3) $M'$ descends to an
isomorphism
$$
M: G\times V@>{(id_G,m)}>> G\times V.
$$
The fact that $m'$ is a $G$-action is expressed by certain commutative
diagrams, which can be made to be diagrams of schemes \'etale over
$G\times G\times X'$.  By (5.3) these diagrams descend to $X$, showing that $V$
is a $G$-scheme and $V'\to V$ is a $G$-morphism.

Finally, pick a point $v'\in V$ and let $v\in V$, $x'\in X'$ and $x\in X$
denote its images. Since $V'\to V$ and $X'\to X$ are universal
homeomorphisms, we see that
$$
\align
&\red\stab(v')=\red\stab(v)\qtq{and} \red\stab(x')=\red\stab(x), \qtq{hence}\\
&\red (\stab(x')/\stab(v'))=\red(\stab(x)/\stab(v)).
\endalign
$$
The fiber of $V'\to X'$ (resp. of $V\to X$) over $x'$ (resp over $x$)
is $\stab(x')/\stab(v')$ (resp. $\stab(x)/\stab(v)$). Both of these
morphisms are \'etale, thus the fibers are smooth. They have the same reduced
scheme structure, hence $\stab(x')/\stab(v')=\stab(x)/\stab(v)$.  This
implies that
$$
\stab(v')= \stab(x')\qtq{iff} \stab(v)=\stab(x).
$$
This completes the proof of (7.8).\qed\enddemo

\head 8.  Push-out Diagrams
\endhead

In this section we  study push-out diagrams for algebraic spaces. These
diagrams are used in going from normal to weakly-normal spaces.

\demo{8.1 Definition of Push-out Diagrams}  The general problem is the
following. Given a diagram
$$
\CD
X@>g_1>> Y_1\\
@Vg_2VV @.\\
Y_2 @.
\endCD\tag 8.1.1
$$
we would like to complete it to a  commutative diagram
$$
\CD
X@>g_1>> Y_1\\
@Vg_2VV @VVs_1V\\
Y_2 @>s_2>> Z.
\endCD\tag 8.1.2
$$

In the category of algebraic spaces over $S$ this is always possible, take for
instance $Z=S$.
We are interested in finding an algebraic space $Z$ which is ``large".

In all the cases that we consider,
$g_1$  and $g_2$ are  finite, and  we would like $s_1$  and $s_2$ to be finite
as well. In such a case we say that the push-out diagram has {\it finite
morphisms}.

Assume that we have a diagram (8.1.2) with finite
morphisms. Then there is a {\it universal
push-out}, given by
$$
Z^u\deq\spec_Z\ker[\o_{Y_1}+\o_{Y_2}@>{g_1^*-g_2^*}>> \o_X].\tag 8.1.3
$$
This implies that the construction of   universal
push-out diagrams is a local question over $Z$:
\enddemo

\proclaim{8.2 Lemma}   Assume that (8.1.2) is a universal push-out
diagram with finite morphisms. Let $Z'\to Z$ be a flat morphism. Then the
diagram obtained by base change
$$
\CD
X\times_ZZ'@>g'_1>> Y_1\times_ZZ'\\
@Vg'_2VV @VVs'_1V\\
Y_2\times_ZZ' @>s'_2>> Z'
\endCD
$$
is again a universal push-out diagram
with finite morphisms.

Thus the construction of universal
push-out diagrams with finite morphisms is a local question in the \'etale
topology of
$Z$.\qed
\endproclaim

Unfortunately, in may cases finite push-out diagrams do not exist.
 The following examples give a sample of   push-out diagrams
with bad properties.

\demo{8.3 Examples} Let $S=\spec k$ for a field $k$.

(8.3.1) Consider the diagram
$$
\CD
\a^1@>{z^2}>> \a^1\\
@V{z^2+z}VV @.\\
\a^1 @.
\endCD
$$
The only possible push-out is $Z=\spec k$, unless $\chr k=2$.

(8.3.2) (cf. \cite{Holmann63, p. 342}) Let $X=\spec k[x,y]$, $Y_1=\spec
k[x,y^2,y^3]$ and
$Y_2=\spec k[x+y,y^2,y^3]$. $g_i$ are the natural morphisms; both are
finite and universal homeomorphisms.

If $\chr k=0$ then the universal push-out is given by
$$
Z=\spec (k[x,y^2,y^3]\cap k[x+y,y^2,y^3])=\spec k[x^iy^2,x^iy^3\:
i=1,2,\dots],
$$
which is not of finite type. The induced morphisms $s_i$ are not finite.

The dependence on the characteristic becomes clear if we note that, as
$k$-vector spaces
$$
\align
k[x,y^2,y^3]=&k[x^iy^2,x^iy^3\: i=1,2,\dots]+\sum_{i>0} kx^i,\qtq{and}\\
k[x+y,y^2,y^3]=&k[x^iy^2,x^iy^3\: i=1,2,\dots]+\sum_{i>0}
k(x^i+ix^{i-1}y).\\
\endalign
$$
The second summands are disjoint in characteristic zero but have large
intersection in positive characteristic.

(8.3.3) In (8.3.2) we may replace the polynomial rings with the
corresponding power series rings. This shows that even in the category of
algebraic or analytic spaces (in characteristic zero), there is no push-out
with
$s_i$ finite.
\enddemo

In positive characteristic the situation is much better:

\proclaim{8.4 Lemma} In the push-out diagram (8.1.1) assume that
all algebraic spaces are of finite type over a scheme of positive
characteristic. Assume furthermore that
$g_1$ is a finite and universal homeomorphism. (No assumption on $g_2$.)  Then
there is a push-out  diagram
$$
\CD
X@>g_1>> Y_1\\
@Vg_2VV @VVs_1V\\
Y_2 @>s_2>> Z,
\endCD
$$
where $s_2$ is also a finite and universal homeomorphism.
\endproclaim

\demop The proof relies on properties of the Frobenius morphism (6.2). Choose
$q$ such that
$F^q$ factors as
$$
F^q:X@>g_1>> Y_1@>\bar g_1>> X^{(q)}.
$$
Set $Z=Y_2^{(q)}$, $s_2\deq F^q:Y_2\to Y_2^{(q)}$ and $s_1\deq g_2^{(q)}\circ
\bar g_1$.\qed\enddemo

\demo{8.5 Example} In the situation of (8.4) the universal push-out may be
rather complicated if $g_2$  is not finite. Consider for instance the
  diagram
$$
\CD
\spec k[x]@>>> \spec k[x,y]/(y^2)\\
@VVV @.\\
\spec k @. @..
\endCD
$$
The construction of (8.4) gives $Z=\spec k$. The uiversal push-out exists
and is given by $Z^u=\spec_k k[y,xy,x^2y,\dots]/(y^2)$. $Z^u$ is a
non-Noetherian scheme.
\enddemo

The most useful general result about push-out diagrams with  finite morphisms
is the following. As usual, we assume that the base scheme is excellent.
The original more restrictive assumptions  in \cite{Artin70} can be relaxed
in view of  recent progress in the approximation property of Henselian rings
\cite{Popescu85; Andr\'e94; Spivakovsky94}.

\proclaim{8.6 Theorem} \cite{Artin70, 6.1} In the push-out diagram (8.1.1)
assume  that
$g_1$ is a closed immersion and $g_2$ is finite. Then there is a push-out
diagram
$$
\CD
X@>g_1>> Y_1\\
@Vg_2VV @VVs_1V\\
Y_2 @>s_2>> Z,
\endCD
$$
where $s_2$ is a closed immersion and $s_1$ is finite.\qed
\endproclaim

As (8.3.2) shows, the analog of (8.4) fails in characteristic zero.
In order to get a general push-out result, we concentrate
on finite and  universal homeomorphisms which are isomorphisms
in characteristic  zero.  First we give a  ring theoretic
characterization:

\proclaim{8.7 Lemma} Let $B$ be a Noetherian ring
and $B_1\subset B$ a subring such that $B$ is finite over $B_1$. The following
are equivalent:

(8.7.1)
$\spec B\to
\spec B_1$ is a finite and  universal homeomorphism which is an isomorphism
at all points where the characteristic of the residue field is zero.

(8.7.2)
There is an $m>0$ such that
$b^m\in B_1$ for every $b\in B$.
\endproclaim

\demop   Assume (8.7.2) and let $P<B_1$ be any prime ideal. Then
  $ B/PB$ is an Artin algebra over  $B_1/P$  such that if $b\in B/PB$
then $b^m\in B_1/P$. This implies that $ B/PB$ is local  with maximal
ideal $\root\of{PB}$ and
$B_1/P\subset B/\root\of{PB}$ is a purely inseparable field extension. Thus by
(5.2),
$\spec B\to
\spec B_1$ is a finite and  universal homeomorphism which is an isomorphism
at all points where the characteristic of the residue field is zero.

Conversely, assume (8.7.1).
$B_1$ is also Noetherian by (3.14.1) and $B/B_1$ is  a finite
$B_1$-module. By assumption all associated primes of $B/B_1$ have positive
residue characteristic, thus there is an integer $N$ such that $NB\subset
B_1$. We can replace the inclusion $B_1\subset B$
by the inclusion $B_1/NB\subset B/NB$. Thus it is sufficient to prove (8.7) in
case all the rings are $\z/(N)$-algebras.

$B$ is a direct sum of $\z/(p_i^{n_i})$-algebras $B^i$ where
$N=\prod p_i^{n_i}$.  Correspondingly, $B_1$ decomposes as a sum of
subalgebras $B^i_1\subset B^i$.

Assume that for each $i$ we find $m_i$ such that
$b^{m_i}\in B^i_1$ for every $b\in B^i$. Then $m=\prod m_i$
works for $B$. Thus we are further reduced to the case when
all the rings are $\z/(p^n)$-algebras.

First consider   $B_1/(B_1\cap pB)\subset B/pB$. The induced map on
the spectra is a finite and
universal homeomorphism, hence by (6.6) there is an $r$ such that
$b^{p^r}\in B_1+pB$ for every $b\in B$.

Let $b_1+pb_2$ be a typical element of $B_1+pB$. Then
$$
(b_1+pb_2)^{p^n}=b_1^{p^n}+\sum_{i>0} p^i\binom{p^n}{i}b_1^{p^n-i}b_2^i.
$$
It is easy to see that $p^n$ divides $p^i\binom{p^n}{i}$ for $i>0$, thus
$(b_1+pb_2)^{p^n}=b_1^{p^n}\in B_1$.\qed\enddemo

\proclaim{8.8 Corollary} (8.8.1)  Assume that (8.1.2) is a universal
push-out diagram with finite morphisms such that  $g_1$ is a finite and
universal homeomorphism which is an isomorphism
at all points where the characteristic of the residue field is zero. Then $s_2$
is also a finite and universal
homeomorphism which is an isomorphism
at all points where the characteristic of the residue field is zero.

(8.8.2)  Assume that (8.1.1) is a  push-out
diagram with finite morphisms such that  $g_1$ is a finite and universal
homeomorphism which is an isomorphism
at all points where the characteristic of the residue field is zero.  Then the
construction of the universal push-out diagram is a local problem in the
\'etale topology of $Y_2$.
\endproclaim

\demop  By (8.2) it is sufficient to consider (8.8.1) in  case  $Z$ is an
affine scheme.  Then $X$ and the $Y_i$ are also affine. By (8.7) there is an
$m>0$ such that
$$
b\in H^0(X,\o_X)\Rightarrow b^m\in \im\left[H^0(Y_1,\o_{Y_1})@>g_1^*>>
H^0(X,\o_X)\right].
$$
Thus, by the construction (8.1.3) of the universal push-out,
$$
b\in H^0(Y_2,\o_{Y_2})\Rightarrow b^m\in \im\left[H^0(Z,\o_{Z})@>s_2^*>>
H^0(Y_2,\o_{Y_2})\right].
$$
By (8.7) this implies that $s_2$ is   a finite and universal
homeomorphism which is an isomorphism
at all points where the characteristic of the residue field is zero.

In order to see (8.8.2), let $V\to Y_2$ be an \'etale morphism. Since $s_2$ is
a  finite and universal
homeomorphism, $V$ descends to an \'etale cover $V_Z\to Z$ such that
$V=V_Z\times_ZY_2$ (5.4). Thus (8.2)  and (8.8.1) imply (8.8.2).\qed\enddemo

This enables us to generalize (8.4) to algebraic spaces over an excellent base
scheme
$S$:

\proclaim{8.9 Lemma}  In the push-out diagram
$$
\CD
X@>g_1>> Y_1\\
@Vg_2VV \\
Y_2 @.
\endCD
$$
 assume
that all algebraic spaces are of finite type over $S$. Assume
that
$g_2$ is finite and
$g_1$ is a finite and universal homeomorphism which is an isomorphism
at all points where the characteristic of the residue field is zero.

  Then there is a push-out  diagram
$$
\CD
X@>g_1>> Y_1\\
@Vg_2VV  @VVs_1V\\
Y_2 @>s_2>>  Z,
\endCD
$$
where $s_1$ is finite and
$s_2$ is a finite and universal homeomorphism which is an isomorphism
at all points where the characteristic of the residue field is zero.
\endproclaim

\demop $g_1$ can be factored as
$$
g_1:X@>h_1>> X_1@>g'_1>> Y_1
$$
where $g'_1$ is a closed immersion and  $\o_{X_1}\to (h_1)_*\o_X$ is
injective.  We construct the push-out in two steps:
$$
\CD
X@>h_1>> X_1@>g'_1>> Y_1\\
@Vg_2VV  @VVs'_1V @VVs_1V\\
Y_2 @>s_2>>  W @>s'_2>>  Z.
\endCD
$$

By (8.8.2) the construction of $W$ is  local in the \'etale topology of $Y_2$.
Let $\spec C\subset S$ be an affine open  set.
 Let $U$ be an affine scheme of finite type over $\spec C$ and $U\to Y_2$
an \'etale morphism. Let   $V= U\times_{Y_2}X$. Then $V\to U$ is finite,
hence $V$ is affine. Also, $V\to X$ is
an \'etale morphism. Since
$g_1$ is a universal homeomorphism, $V$ descends to an affine scheme $V_1$
which is
\'etale over
$X_1$ (5.4). Corresponding to the morphisms of affine schemes
$U@<<< V@>>> V_1$ we obtain
morphisms of rings
$$
A@>g_2^*>> B\overset{h_1^*}\to{\hookleftarrow} B_1.
$$
$B_1\subset B$ is a subring such that $\spec B\to
\spec B_1$ is a finite and universal homeomorphism which is an isomorphism
at all points where the characteristic of the residue field is zero.
Choose $m$ as in (8.7.2).  Let $A_2\subset A$ be the $C$-subalgebra generated
by the
$m$-th powers of all elements of $A$.   $A$ is a finitely generated
$C$-algebra, thus $A$ is finite over $A_2$.
By  (8.7), $g_2(A_2)\subset B_1$.
Set
$$
A_1\deq \{a\in A\vert g_2(a)\in B_1\}, \qtq{and} U_1\deq \spec A_1.
$$
Thus we obtain a commutative diagram
$$
\CD
V@>g_1>> V_1\\
@Vg_2VV  @VVs_1V\\
U @>s_2>>  U_1,
\endCD
$$
which gives the universal push-out. By (8.2.3)  the schemes $U_1$ give an
\'etale covering of $W$.

Once we have $W$, we can apply (8.6) to the diagram
$$
\CD
 X_1@>g'_1>> Y_1\\
  @VVs'_1V @.\\
 W @. @.
\endCD
$$
to obtain $Z$.
 \qed\enddemo

In some applications (8.9) is used in the following slightly more general form:

\proclaim{8.10 Lemma}  In the push-out diagram
$$
\CD
X@>g_1>> X_1@>h_1>> Y_1\\
@Vg_2VV @. @.\\
Y_2 @.@.
\endCD
$$
 assume
that  $g_2$
is finite and $h_1$ is a closed immersion. Assume furthermore that
$g_1$ is a finite and universal homeomorphism which is an isomorphism
at all points where the characteristic of the residue field is zero.

  Then there is a push-out  diagram
$$
\CD
X@>g_1>> X_1@>h_1>> Y_1\\
@Vg_2VV @. @VVs_1V\\
Y_2 @>>> {} @>>> Z,
\endCD
$$
where all morphisms are finite.
\endproclaim

\demop  By (8.9) there is  a push-out  diagram
$$
\CD
X@>g_1>> X_1\\
@Vg_2VV  @VVh_2V\\
Y_2 @>>> Y'_2,
\endCD
$$
where all morphisms are finite. Then we can apply (8.6) to the diagram
$$
\CD
X_1@>h_1>> Y_1\\
@Vh_2VV @. \\
Y'_2 @.
\endCD
$$
to find $Z$ as required. \qed\enddemo

\bigpagebreak

Department of Mathematics,

University of Utah,

Salt Lake City
UT 84112, USA

kollar\@{}math.utah.edu

\end